\documentclass[preprint]{elsarticle}

\usepackage{hyperref}

\usepackage{lineno}
\modulolinenumbers[5]

\usepackage{amssymb,latexsym,amsmath,amsthm,graphicx,epsfig,xcolor}
\usepackage{mathtools}

\usepackage{mathabx}
\usepackage{multirow}

\newtheorem{theorem}{Theorem}[section]
\newtheorem{lemma}[theorem]{Lemma}

\newtheorem{corollary}[theorem]{Corollary}
\newtheorem{definition}[theorem]{Definition}

\newcommand{\ip}[2]{\langle#1,#2\rangle}
\newcommand{\abs}[1]{| #1 |}
\newcommand{\norm}[1]{\lVert#1\rVert}

\definecolor{schrift}{rgb}{0,0,.4} 

\newcounter{cnt}
\newcommand{\addPDFfig}[1]{
	\setcounter{cnt}{#1}
	\stepcounter{cnt}
	
	\includegraphics[width=.99\linewidth, page=\value{cnt}]{figuresRevised.pdf}
}

\def\supp{{\text{\rm supp }\!}}

\def\E{\mathcal{E}}
\def\F{\mathcal{F}}

\def\N{\mathbb{N}}
\def\M{\mathcal{M}}

\def\Z{\mathbb{Z}}
\def\Q{\mathcal{Q}}
\def\R{\mathbb{R}}

\def\P{\mathcal{P}}

\def\X{\mathcal{X}}

\def\E{\mathcal{E}}

\def\f{\frac}
\allowdisplaybreaks

\newcommand{\cc}{*_{\text{3d}}}
\renewcommand{\vec}[1]{\boldsymbol{#1}} 
\newcommand{\fvec}{\vec{f}}
\newcommand{\m}{\vec{m}}
\newcommand{\y}{\vec{y}}
\newcommand{\rvec}{\vec{r}}
\newcommand{\etab}{\vec{\eta}}
\newcommand{\uShc}{\vec{u}_{j,\ell,d}}
\newcommand{\Wvec}{\vec{W}} 
\newcommand{\Vvec}{\vec{V}} 

\DeclareMathOperator{\RadonD}{\vec{\mathcal{R}}}

\DeclareMathOperator{\CylShTD}{\vec{\mathcal{S}}}
\newcommand{\proj}[1]{\text{proj}_{#1}} 

\begin{document}

\begin{frontmatter}
	
	\title{Efficient representation of spatio-temporal data using cylindrical shearlets\tnoteref{journal}}
	
	\tnotetext[journal]{Published in the Journal of Computational and Applied mathematics 429 (2023) 115206: \href{https://doi.org/10.1016/j.cam.2023.115206}{doi:10.1016/j.cam.2023.115206}.}
	
	\author[TBaddress]{Tatiana A.~Bubba}

    \author[GEaddress]{Glenn Easley}
	
	\author[THaddress]{Tommi Heikkil\"a}
	
	\author[DLaddress]{Demetrio Labate\corref{mycorrespondingauthor}}
	\cortext[mycorrespondingauthor]{Corresponding author}
    \ead{dlabate@math.uh.edu}
	
	\author[JRaddress]{Jose P.~Rodriguez Ayllon}

	\address[TBaddress]{Department of Mathematical Sciences, University of Bath, Claverton Down, Bath \mbox{BA2 7AY}, United Kingdom}
		\address[GEaddress]{Applied Physics Laboratory, Johns Hopkins University, 11100 Johns Hopkins Road
		Laurel, Maryland 20723, USA }
	\address[THaddress]{Department of Mathematics and Statistics, University of Helsinki, Pietari Kalmin katu 5, 00014 Helsinki, Finland}

	\address[DLaddress]{Department of Mathematics, University of Houston, 651 Phillip G Hoffman
		Houston, Texas 77204-3008, USA }

	\address[JRaddress]{Universidad Mayor de San Andres, Av. Villazon 1995 Monoblock Central, La Paz, Bolivia}

	\begin{abstract}
		Efficient representations of multivariate functions are critical for the design of state-of-the-art methods of data restoration and image reconstruction. In this work, we consider the representation of spatio-temporal data such as temporal sequences (videos) of 2- and 3-dimensional images,  where conventional separable representations are usually very inefficient, due to their limitations in handling the geometry of the data. To address this challenge, we define a class $\E(A) \subset L^2(\R^4)$ of functions of 4 variables dominated by hypersurface singularities in the first three coordinates that we apply to model 4-dimensional data corresponding to temporal sequences (videos) of 3-dimensional objects.
		
		To provide an efficient representation for this type of data, we introduce a new multiscale directional system of functions based on cylindrical shearlets and prove that this new approach achieves superior approximation properties with respect to conventional multiscale representations. We illustrate the advantages of our approach by applying a discrete implementation of the new representation to a challenging problem from dynamic tomography. Numerical results confirm the potential of our novel approach with respect to conventional multiscale methods.
	\end{abstract}
	
	\begin{keyword}
		 dynamic tomography \sep multiscale analysis  \sep shearlets \sep spatio-temporal data \sep sparse approximations \sep regularization
		 \MSC[2010]  42C10 \sep
		 42C40  \sep
		 92C55 \sep 47A52
	\end{keyword}
	
\end{frontmatter}


\section{Introduction}

Sparse representations of multivariate functions have been remarkably successful in applied mathematics and signal processing,  with applications ranging from image denoising and inpainting through medical image reconstruction and feature extraction being proposed during the last decade. A multiplicity of such constructions were introduced  to deal with different types of multidimensional data and signal processing tasks, including curvelets \cite{CD2004}, shearlets \cite{GL1,kutyniok2011}, bandlets \cite{mallat2009}, scattering wavelets \cite{bruna2013}, bendlets \cite{LPS2019}, parabolic molecules \cite{grohs2014}, directional multivariate wavelets \cite{bergmann2015} and directional framelets \cite{han2019}. The key observation underpinning such constructions is that representations of multivariate functions that  capture the fundamental geometry of data result in superior approximation properties that can be translated into improved algorithms for signal processing applications. For instance, shearlets, which are defined as  well-localized anisotropic waveforms  ranging over multiple scales, location and orientations in $L^2(\R^2)$, are especially designed to represent edge discontinuities. As a result, they provide optimally sparse approximations, in a precise sense, for cartoon-like images - a class of piecewise smooth functions that is used to model a large class of natural images - outperforming conventional multiscale representations. Such approximation properties were critically exploited to develop successful numerical algorithms for signal processing and medical imaging~\cite{bubba2019,colonna2010,easley2008sparse,king2014,yi2009}.

In this paper, we introduce a new construction of {\it cylindrical shearlets} on $L^2(\R^4)$ aimed at the efficient representation of spatio-temporal data, that is, temporal sequences (or videos) of 3-dimensional objects. Our approach is especially motivated by dynamic computed tomography (CT), a medical imaging technique whose goal is to reconstruct 3-dimensional image sequences where the main focus is the dynamic of the living human body \cite{Bonnet2003} for applications such as cardiac imaging or image-guided interventional medical procedures. A main challenge in dynamic CT reconstruction is that, due to technical or physical constraints, data are often heavily undersampled causing the inverse problem associated with the reconstruction task to be potentially ill-posed. For instance, many dynamic CT scenarios involve the use of contrast tracers and full X-ray scans are too slow to capture the movement of the tracer (e.g., iodine) in the imaging windows. The most common remedy for reducing the duration of the imaging as well as the radiation dose consists in lowering the number of scanning angles leading to an undersampled reconstruction problem. 

Some of the authors of this paper have recently shown that one can successfully address the undersampled reconstruction problem in dynamic CT by taking advantage of appropriate sparse data representations~\cite{bubba2020sparse}. In particular, 3-dimensional shearlets were successfully applied to develop an improved algorithm for the reconstruction of 2-dimensional time frames in sparse dynamic tomography by exploiting their superior approximation properties of 3-dimensional data. However, their method does not apply directly to the `full' dynamic CT problem of reconstructing 3-dimensional time frames, that is, 4-dimensional data. To deal with such task, here we introduce cylindrical shearlets on $L^2(\R^4)$ as a collection of well-localized waveforms ranging over multiple scales, locations and orientations on $\R^4$. To better adapt the geometry of this representation to the characteristics of spatio-temporal data - under the simplifying assumptions that such data are dominated by hyper-surface discontinuities in the three spatial coordinates - we will assume that our representation has directional sensitivity with respect to the 3 spatial coordinates but not along the time coordinate.

Our main theoretical result in this paper is that this new construction provides highly sparse representations for the class of 4-dimensional cylindrical cartoon-like functions - the simplified model we adopt for spatio-temporal data - outperforming more conventional representations. Next, to illustrate the potential of our new construction in numerical applications, we consider a problem of undersampled reconstruction in dynamic tomography using synthetic data. Our numerical results show that our algorithm for dynamic CT reconstruction based on cylindrical shearlets improves the reconstruction quality as compared to similar methods based on conventional wavelets.

Finally,  we remark that shearlets have been already applied in (2d+1) video applications. For instance, conventional 3d-shearlets were employed in~\cite{NL2012} to provide efficient video representations and, more recently, to detect relevant space-time features of videos in~\cite{Malafronte2018}. However, the cylindrical shearlets we consider in this paper are derived from a very different construction that handles spatial and temporal coordinates with different geometric sensitivities. As already indicated by some of the authors in~\cite{LPGE00} and further argued in this paper, this construction entails distinct mathematical properties with respect to conventional shearlets and significant potential advantages in the context of spatio-temporal data.

\subsection{Sparse 4-dimensional representations}
To explain the significance of our new representation, we start with  a heuristic argument showing why cylindrical shearlets are expected to be especially effective in representing a compactly supported piecewise regular function $f$ of four variables with discontinuities in the first 3 spatial coordinates. To keep this explanation at an intuitive level, we will postpone the precise definition of the class $\E(A)$ of cylindrical cartoon-like functions to Sec.~\ref{s.main}.

We start by examining the {\it 4d wavelet expansion} of $f$ using a Parseval frame of wavelets $\{ \varphi_{j,k} (x) = 2^{4j} \varphi(2^{2j} x - k) :j \in \Z, k \in \Z^4 \} \subset L^2(\R^4)$ where $\varphi$ is well localized. We choose $2^{2j}$ as dilation factor rather than $2^j$  to be consistent with the cylindrical shearlet representation. An element $\varphi_{j,k}$ of the wavelet system at scale of $2^{-2j}$ is essentially supported on a box of size $2^{-2j} \times 2^{-2j} \times 2^{-2j} \times 2^{-2j}$. Since the surface of discontinuity of $f$ has finite volume in the 4d-space, there are approximately $2^{6j}$ wavelet coefficients $F_{j,k} (f)=\ip{f}{ {\varphi_{j,k}}}$ associated with this surface, while the remaining coefficients are negligible at fine scales. A direct computation shows that
$$\int_{\R^4}\abs{\varphi_{j,k}(x)} \, dx=2^{4j}\int_{\R^4}\abs{\varphi(2^{2j} x - k)} \, dx=2^{-4j}\int_{\R^4}\abs{\varphi(y)} \, dy \le c \, 2^{-4j},$$
for a constant $c>0$. Hence, at scale of $2^{-2j}$, we have 
$$\abs{F_{j,k} (f )} \leq \norm{f }_{\infty}\norm{\varphi_{j,k}}_{L^1} \leq c \,2^{-4j}.$$
For brevity, here and in the following we use the convention that the same letter $c$ or $C$ may denote different uniform constants.
Thus, letting $N = 2^{6j}$, the $N$-th largest wavelet coefficient of $f$ in magnitude, denoted by $|F(f)|_N$, is  bounded by $O(N^{-\frac{2}{3}} )$. Hence, if $f^{(wav)}_N$
is the approximation of $f$ obtained by taking the $N$ largest coefficients (in absolute value)  of its wavelet expansion, we have
$$\norm{f-f^{(wav)}_N}^2_{L^2}\leq \sum_{\mu>N}\abs{F(f)}_\mu^2\leq c \,  N^{-\frac{1}{3}}.$$

Next, we examine the {\it  4d cylindrical shearlet expansion} of $f$. The elements of the cylindrical shearlet system $\psi_{ j,\ell,k}$ are essentially of the form $2^{3j} \psi(A^jB_{\ell}x-k)$, where $\psi$ is a bounded well-localized function, $A$ is a diagonal matrix with factors $(2^{2j}, 2^j, 2^j, 2^{2j})$ and $B_{\ell}$ is an appropriate shear matrix. As a result, a direct estimate shows that there is a constant $c>0$ such that
$$\int_{\R^4}\abs{\psi_{j,k,\ell}(x)} \, dx=2^{-3j}\int_{\R^4}\abs{\psi(y)} \, dy \, \le c \, 2^{-3j}.$$
Hence, at scale $2^{-2j}$, the cylindrical shearlet coefficients $s_{j,k,\ell}=\ip{f }{\psi_{ j,\ell,k}}$ are bounded by
$$\abs{s_{j,k,\ell} (f )} \leq \norm{f }_{\infty}\norm{\psi_{j,k,\ell}}_{L^1} \leq c \, 2^{-3j}.$$
Each element $\psi_{ j,\ell,k}$ is essentially supported on a parallelepiped of size $2^{-2j}\times 2^{-j}\times 2^{-j} \times 2^{-2j}$ with various orientations controlled by $\ell$ and, due to directional sensitivity and elongated support, the only significant cylindrical shearlet coefficients occur when an element $\psi_{ j,\ell,k}$ is tangent to the surface of discontinuity of $f$. Only about $2^{3j}$ cylindrical shearlet coefficients  are significant. 
Thus, letting $N = 2^{3j}$, the N-th largest cylindrical shearlet coefficient in absolute value, denoted as $|s(f)|_N$, is  bounded by $O(N^{-1})$. If we denote as $f^{(csh)}_N$ the approximation of $f$ obtained by taking the $N$ largest coefficients (in absolute value) of its shearlet expansion, we have
$$\norm{f-f^{(csh)}_N}^2_{L^2}\leq \sum_{\mu>N}\abs{s(f)}_\mu^2\leq c \, N^{-1}.$$
We will show below using a rigorous argument that the estimate above is essentially correct.

\subsection{Outline} The rest of the paper is organized as follows.
In Sec.~\ref{s.sh}, we introduce a new construction of 4-dimensional cylindrical shearlets by generalizing the 3-dimensional construction in~\cite{LPGE00}. In Sec.~\ref{s.main}, we present our sparse approximation results using 4-dimensional cylindrical shearlets, whose proofs are postponed to Appendix~\ref{s.proof}. We present a numerical implementation of 4-dimensional cylindrical shearlets in Sec.~\ref{s.numImplementation} and apply this representation to a problem from dynamic tomography in Sec.~\ref{s.numDynamicCT}.
We finally provide concluding remarks in Sec.~\ref{s.conclusion}.

\section{Cylindrical shearlets}  \label{s.sh}

Cylindrical shearlets were recently introduced by some of the authors~\cite{LPGE00} as a variant of the shearlet construction in the 3-dimensional setting. As remarked above, this construction is motivated by applications where data are dominated by discontinuities that occur perpendicularly to one of the coordinate axes so that it is useful to employ representations that are direction-sensitive with respect to one hyperspace.

In the 4-dimensional setting, we associate cylindrical shearlets to three cylindrical hyper-pyramids  defined as:
\begin{eqnarray*}
 \mathcal{P}_1 &=& \{(\xi_1,\xi_2,\xi_3,\xi_4) \in \R^4: \abs{\tfrac{\xi_2}{\xi_1}} \leq 1, \abs{\tfrac{\xi_3}{\xi_1}} \leq 1   \},\\
  \mathcal{P}_2&=&\{(\xi_1,\xi_2,\xi_3,\xi_4) \in \R^4: \abs{\tfrac{\xi_1}{\xi_2}} \leq 1, \abs{\tfrac{\xi_3}{\xi_2}} \leq 1   \},\\
    \mathcal{P}_3&=&\{(\xi_1,\xi_2,\xi_3,\xi_4) \in \R^4: \abs{\tfrac{\xi_1}{\xi_3}} \leq 1, \abs{\tfrac{\xi_2}{\xi_3}} \leq 1   \}.
  \end{eqnarray*}
\begin{definition}\label{def:cylindrical}
For $d=1,2,3$, a {\it pyramid-based cylindrical shearlet system} associated with the pyramid ${\mathcal{P}}_d$ is a collection of functions 
	\begin{equation} \label{sh.one}
	\{\psi^{(d)}_{j,\ell,k}: j \geq 0,\, \ell = (\ell_1, \ell_2) \in \Z^2, |\ell_1|,|\ell_2|\leq2^j,\, k \in \Z^4 \},
	\end{equation}
where the elements of the system~\eqref{sh.one} are given in the Fourier domain as
\begin{equation}\label{sh.one.f}
\hat{\psi}^{(d)}_{j,\ell,k}(\xi)= \abs{\det A_{(d)}}^{-\frac{j}{2}} W(2^{-2j}\xi) \, V_{(d)}(\xi A_{(d)}^{-j}B_{(d)}^{[-\ell]}) \, e^{2\pi i\xi A_{(d)}^{-j}B_{(d)}^{[-\ell]}k }
\end{equation}
with functions $W, V_{(d)}: \R^4 \to  {\mathbb C}$ to be defined below and the matrices $A_{(d)}$ and $B_{(d)}^{[\ell]}$ given by 
 \newline
	$A_{(1)}=\begin{pmatrix}
	4 & 0 & 0 &0\\
	0 & 2 & 0 &0\\
	0 & 0 & 2 &0 \\
	0& 0& 0 & 4
	\end{pmatrix}$, 
	$A_{(2)}=\begin{pmatrix}
	2 & 0 & 0 & 0\\
	0 & 4 & 0 & 0\\
	0 & 0 & 2 & 0 \\
		0 & 0 & 0 & 4 
	\end{pmatrix}$,
		$A_{(3)}=\begin{pmatrix}
	2 & 0 & 0 & 0\\
	0 & 2 & 0 & 0\\
	0 & 0 & 4 & 0 \\
	0 & 0 & 0 & 4 
	\end{pmatrix}$, \newline
	$B_{(1)}^{[\ell]}=\begin{pmatrix}
	1 & \ell_1 & \ell_2 &0\\
	0 & 1 & 0 &0\\
	0 & 0 & 1 & 0 \\
	0 & 0 & 0 & 1
	\end{pmatrix}$,
	 $B_{(2)}^{[\ell]}=\begin{pmatrix}
	1 & 0 & 0 &0\\
	\ell_1 & 1 & \ell_2 &0\\
	0 & 0 & 1  & 0 \\
	0 & 0 & 0 & 1
	\end{pmatrix}$,
	$B_{(3)}^{[\ell]}
	=\begin{pmatrix}
	1 & 0 & 0 &0\\
	0 & 1 & 0 &0\\
	\ell_1 & \ell_2 & 1  & 0 \\
	0 & 0 & 0 & 1
	\end{pmatrix}$.	
\end{definition}
As we show below, we can choose the functions $W$ and $V_{(d)}$ so that the corresponding system~\eqref{sh.one} is a smooth Parseval frame of $L^2 ({\mathcal{P}}_d\setminus C_0)^\vee$, for $d=1,2,3$ where $C_0 =[-\frac 18, \frac 18]^4$, that is,
$$ f = \sum_{j \ge 0} \sum_{|\ell_1|, |\ell_2| \le 2^j} \sum_{k \in \Z^4} \ip{f}{\psi^{(d)}_{j,\ell,k}} \psi^{(d)}_{j,\ell,k}, $$
for all $f$ in $L^2$ whose Fourier support is contained in ${\mathcal{P}}_d\setminus C_0$; convergence is understood in the $L^2$ norm.

\subsection{Smooth Parseval frame of cylindrical shearlets on $L^2(\R^4)$} \label{ss.parsevalConstruction}
Our construction below extends the original 3-dimensional cylindrical shearlet construction~\cite{LPGE00} by adapting some ideas from the standard shearlet construction~\cite{GL_MMNP}.  

We let $\phi \in L^2(\R)$ be such that $\hat \phi \in C_c^\infty$ with $0 \le \hat \phi \le 1$ and
\begin{equation}\label{eq.phi2}
\hat \phi(u) =1 \text{ if } u \in  [-\frac1{16}, \frac 1{16}], \quad \hat \phi(u) =0  \text{ if } u \in \R \setminus [-\frac18, \frac 18].
\end{equation}
For $\xi = (\xi_1,\xi_2,\xi_3,\xi_4) \in \R^4,$ we let $\widehat{\Phi}(\xi_1,\xi_2,\xi_3,\xi_4)=\widehat{\phi}(\xi_1)\widehat{\phi}(\xi_2)\widehat{\phi}(\xi_3)\widehat{\phi}(\xi_4)$ and we define the window function (in $L^2(\R^4$)
\begin{equation*} 
 W(\xi) = \sqrt{ \widehat \Phi^2(2^{-2}\xi)- \widehat \Phi^2(\xi) }. 
 \end{equation*}
It follows that
\begin{equation*} 
\widehat {\Phi}^2(\xi)+ \sum_{j \ge 0} W^2(2^{-2j}\xi) =1 \, \text{ for } \xi \in \R^4.
\end{equation*} We notice that the functions $W^2_j= W^2(2^{-2j} \cdot)$ are supported in the Cartesian coronae
$$ C_j = [- 2^{2j-1}, 2^{2j-1}]^4 \setminus [- 2^{2j-4}, 2^{2j-4}]^4  \subset \R^4  $$
and that, by adding them up for $j \ge 0$, we obtain a smooth tiling of the frequency space $\R^4$ away from the origin:
\begin{equation*}
\sum_{j \ge 0} W^2(2^{-2j}\xi) =1 \,\,\, \text{ for } \, \xi  \in \R^4 \setminus [-\tfrac 18, \tfrac 18]^4.
\end{equation*}
In addition, we let $v\in C^\infty( \R)$ be such that supp$(v) \subset [-1,1]$ 
\begin{equation*} 
|v(u-1)|^2 +|v(u)|^2 +|v(u+1)|^2 =1  \quad \text{for } |u| \leq 1.
\end{equation*}
It is shown in~\cite{GL1} that there exist examples of functions $\phi$ and  $v$ satisfying the properties described above.

For $d=1$, observing that $\abs{\det A_{(1)}}=2^{{6}}$, that
$$(\xi_1,\xi_2,\xi_3,\xi_4)A_{(1)}^{-j}B_{(1)}^{[-\ell]}=(2^{-2j}\xi_1,-2^{-2j}\ell_1\xi_1+2^{-j}\xi_2,-2^{-2j}\ell_2\xi_1+2^{-j}\xi_3,2^{-2j}\xi_4)$$ 
 and setting $V_{(1)} = v(\tfrac{\xi_2}{\xi_1}) v(\tfrac{\xi_3}{\xi_1})$, an element of the  system~\eqref{sh.one.f} can be written as
\begin{equation}\label{shear1}
\hat{\psi}^{(1)}_{j,\ell,k}(\xi)= 2^{-3j} W(2^{-2j}\xi) \, v(2^j\tfrac{\xi_2}{\xi_1}-\ell_1) \, 
v(2^j\tfrac{\xi_3}{\xi_1}-\ell_2) \, e^{2\pi i \xi A_{(1)}^{-j}B_{(1)}^{{[-\ell]}}k },
\end{equation}
showing that the Fourier support of $\psi^{(1)}_{j,\ell,k}$ is contained inside the region
\begin{eqnarray}  \label{U.j.ell}
U_{j,\ell} &=& 
\{ \xi \in [-2^{2j-1},2^{2j-1}]^4\setminus [-2^{2j-4},2^{2j-4}]^4: \,\abs{\tfrac{\xi_2}{\xi_1}-\ell_1 2^{-j}} \leq 2^{-j},  
\nonumber \\
&&  \abs{\tfrac{\xi_3}{\xi_1}-\ell_2 2^{-j}}\leq 2^{-j}\} \subset \R^4.
\end{eqnarray}

Similar to conventional 3-dimensional shearlets \cite{GL_MMNP}, we obtain a smooth Parseval frame of cylindrical shearlets for $L^2(\R^4)$ using an appropriate combinations of the pyramid-based systems~\eqref{sh.one}
together with an additional coarse scale system. To ensure that all elements of this combined system are smooth and compactly supported in Fourier domain, we appropriately modify the elements of the shearlet system overlapping the boundaries of the 
regions $\mathcal{P}_1$, $\mathcal{P}_2$ and $\mathcal{P}_3$. Hence a {\it cylindrical shearlet system for $L^2(\R^4)$} is  given by
\begin{eqnarray}
\Psi &=&  \left\{\!\widetilde \psi_{-1,k}\!: \!k \in \Z^4\!\right\} \! \bigcup \! \left\{\!\widetilde{\psi}_{j,\ell,k,d}\!: \! j \ge 0, |\ell_1| \leq 2^j, |\ell_2| < 2^j, k \in \Z^4, d=1,2,3\right\} \nonumber \\
  && \bigcup \!
  \left\{\!\widetilde{\psi}_{j,\ell,k}\!: \! j \ge 0, \ell_1 =\pm  2^j, \ell_2 =\pm  2^j, k \in \Z^4 \! \right\}\!, \label{def.nss}
\end{eqnarray}
consisting of:
\begin{itemize}
\item the {\it coarse-scale cylindrical shearlets} $\{\widetilde \psi_{-1,k} = \phi(\cdot -k): k \in \Z^4\}$, where $\phi$ is given by~\eqref{eq.phi2};
\item the {\it interior cylindrical shearlets} $\{\widetilde{\psi}_{j,\ell,k,d} = \psi_{j,\ell,k}^{(d)}: j \ge 0, |\ell_1|< 2^j, |\ell_2|< 2^j, k \in \Z^4, d=1,2,3\}$, with the functions $\psi_{j,\ell,k}^{(d)}$ given by \eqref{sh.one};
\item the {\it boundary cylindrical shearlets} $\{\widetilde{\psi}_{j,\ell,k,d}: \, j \ge 0, \ell_1 =\pm  2^j, |\ell_2|<2^j, k \in \Z^4 ,\,d=1,2,3\}$ and $\{\tilde{\psi}_{j,\ell,k}:\,j\geq0\,,\ell_1,\ell_2=\pm2^j\,,k\in\Z^4\}$, obtained by joining together slightly modified versions of $\psi_{j,\ell,k}^{(d)}$
and  $\psi_{j,\ell,k}^{(d')}$, $d \ne d'$, for $\ell_1,\ell_2=\pm 2^j$, after that they have been restricted in the Fourier domain to their pyramids $\P_d$, $\P_{d'}$,  respectively. Their precise definition is very similar to \cite[Sec.~3.1]{GL_MMNP}.
\end{itemize}

We remark that, by construction, the boundary shearlets $\{\widetilde{\psi}_{j,\ell,k}:j \ge 0, \ell_1,\ell_2 =\pm  2^j, k \in \Z^3\}$ are compactly supported in Fourier domain. In addition, we can show that they are smooth in Fourier domain using essentially the same argument as \cite[Sec.~3.1]{GL_MMNP}.

For $j\geq1$, $\ell=(\ell_1,\ell_2)$, $\ell_1=\pm 2^j$, $|\ell_2|<2^j$ we define
\[ (\tilde{\psi}_{j,\ell,k,1})^{\wedge}(\xi)=\begin{cases} 
      2^{-3j-4}W(2^{-2j}\xi) \, v(2^j\f{\xi_2}{\xi_1}-\ell_1) \, v(2^j\f{\xi_3}{\xi_1}-\ell_2) \, h_1(\xi)
      \quad \text{ if } \xi\in\P_1 \\
      2^{-3j-4}W(2^{-2j}\xi) \, v(2^j\f{\xi_1}{\xi_2}-\ell_1) \,v(2^j\f{\xi_3}{\xi_2}-\ell_2) \,h_1(\xi)
      \quad \text{ if } \xi\in\P_2
   \end{cases}
\]
\[ (\tilde{\psi}_{j,\ell,k,2})^{\wedge}(\xi)=\begin{cases} 
      2^{-3j-4}W(2^{-2j}\xi) \, v(2^j\f{\xi_1}{\xi_2}-\ell_1) \, v(2^j\f{\xi_3}{\xi_2}-\ell_2) \, h_2(\xi)
      \quad \text{ if } \xi\in\P_2 \\
      2^{-3j-4}W(2^{-2j}\xi) \, v(2^j\f{\xi_1}{\xi_3}-\ell_1) \, v(2^j\f{\xi_2}{\xi_3}-\ell_2) \, h_2(\xi) 
      \quad \text{ if } \xi\in\P_3
   \end{cases}
\]
\[ (\tilde{\psi}_{j,\ell,k,3})^{\wedge}(\xi)=\begin{cases} 
      2^{-3j-4}W(2^{-2j}\xi) \, v(2^j\f{\xi_2}{\xi_1}-\ell_2) \, v(2^j\f{\xi_3}{\xi_1}-\ell_1) \, h_3(\xi) 
      \quad \text{ if } \xi\in\P_1 \\
      2^{-3j-4}W(2^{-2j}\xi) \, v(2^j\f{\xi_1}{\xi_3}-\ell_1) \, v(2^j\f{\xi_2}{\xi_3}-\ell_2) \, h_3(\xi)
      \quad \text{ if } \xi\in\P_3
   \end{cases}
\]
where $h_d(\xi) = e^{2\pi i \xi 2^{-2}A_{(d)}^{-j}B_{(d)}^{[-(\ell_1,\ell_2)]}k}$, for $d=1,2,3$.

Similarly, for $j\geq1$, $\ell_1,\ell_2=\pm 2^j$ we define
\[ (\tilde{\psi}_{j,\ell,k})^{\wedge}(\xi)=\begin{cases} 
      2^{-3j-4}W(2^{-2j}\xi)\, v(2^j\f{\xi_2}{\xi_1}-\ell_1) \,v(2^j\f{\xi_3}{\xi_1}-\ell_2) \,h_1(\xi) \quad \text{ if } \xi\in\P_1 \\
      2^{-3j-4}W(2^{-2j}\xi)\, v(2^j\f{\xi_1}{\xi_2}-\ell_1) \,v(2^j\f{\xi_3}{\xi_2}-\ell_2)  \,h_1(\xi) \quad \text{ if } \xi\in\P_2\\
      2^{-3j-4}W(2^{-2j}\xi)\, v(2^j\f{\xi_1}{\xi_3}-\ell_1) \, v(2^j\f{\xi_2}{\xi_3}-\ell_2)  \,h_1(\xi) \quad \text{ if } \xi\in\P_3.
   \end{cases}
\]
For $j=0$, $\ell_1=\pm1$, we define
\[ (\tilde{\psi}_{0,\ell_1,0,k,1})^{\wedge}(\xi)=\begin{cases} 
      W(\xi)\, v(\f{\xi_2}{\xi_1}-\ell_1) \,v(\f{\xi_3}{\xi_1})  \, e^{2\pi i \xi k}\quad \text{ if } \xi\in\P_1 \\
      W(\xi) \, v(\f{\xi_1}{\xi_2}-\ell_1) \, v(\f{\xi_3}{\xi_2}) \, e^{2\pi i \xi k}\quad \text{ if } \xi\in\P_2
   \end{cases}
\]
\[ (\tilde{\psi}_{0,\ell_1,0,k,2})^{\wedge}(\xi)=\begin{cases} 
      W(\xi) \,v(\f{\xi_1}{\xi_2}) \,v(\f{\xi_3}{\xi_2}-\ell_1)  \, e^{2\pi i \xi k}\quad \text{ if } \xi\in\P_2 \\
      W(\xi) \, v(\f{\xi_1}{\xi_3}) \, v(\f{\xi_2}{\xi_3}-\ell_1) \, e^{2\pi i \xi k}\quad \text{ if } \xi\in\P_3
   \end{cases}
\]
\[ (\tilde{\psi}_{0,\ell_1,0,k,3})^{\wedge}(\xi)=\begin{cases} 
      W(\xi) \,v(\f{\xi_2}{\xi_1}) \,v(\f{\xi_3}{\xi_1}-\ell_1) \, e^{2\pi i \xi k}\quad \text{ if } \xi\in\P_1 \\
      W(\xi)\, v(\f{\xi_1}{\xi_3}-\ell_1) \, v(\f{\xi_2}{\xi_3}) \, e^{2\pi i \xi k}\quad \text{ if } \xi\in\P_3.
   \end{cases}
\]

For $j=0$, $\ell_1,\ell_2=\pm1$, we define
\[ (\tilde{\psi}_{0,\ell_1,\ell_2,k})^{\wedge}(\xi)=\begin{cases} 
      W(\xi)\,v(\f{\xi_2}{\xi_1}-\ell_1) \, v(\f{\xi_3}{\xi_1}-\ell_2)\,  e^{2\pi i \xi k}\quad \text{ if } \xi\in\P_1 \\
      W(\xi)\,v(\f{\xi_1}{\xi_2}-\ell_1) \,v(\f{\xi_3}{\xi_2}-\ell_2)\,  e^{2\pi i \xi k}\quad \text{ if } \xi\in\P_2 \\
      W(\xi)\, v(\f{\xi_1}{\xi_3}-\ell_1) \,v(\f{\xi_2}{\xi_3}-\ell_2)\,  e^{2\pi i \xi k}\quad \text{ if } \xi \in\P_3
   \end{cases}
\]

We have the following result whose proof is similar to \cite{GL_MMNP}.
\begin{theorem} 
The  shearlet system $\Psi \subset L^2(\R^4)$, given by~\eqref{def.nss}, is a Parseval frame for $L^2(\R^4)$. Furthermore, the elements of this system are $C^\infty$ and compactly supported in the Fourier domain.
\end{theorem}

For simplicity, in the following we will denote the cylindrical system of shearlets in~\eqref{def.nss}  as
\begin{equation}\label{shearSystem}
   \Psi = \{\tilde{\psi}_\mu\,:\,\mu\in M\},
\end{equation}
where $M=M_C\cup M_I\cup M_B$ are the indices associated with
\begin{itemize}
    \item $M_C=\{(j,k)\,:\,j=-1,\,k\in\Z^4\}$ coarse-scale shearlets,
    \item $M_I=\{(j,\ell_1,\ell_2,k,d)\,:\,  j \ge 0, |\ell_1|< 2^j, |\ell_2|< 2^j, k \in \Z^4, d=1,2,3  \}$ interior shearlets,
    \item $M_B=\{(j,\ell_1,\ell_2,k,d)\,:\,  j \ge 0, |\ell_1|=\pm 2^j, |\ell_2|< 2^j, k \in \Z^4, d=1,2,3 \}\cup\{(j,\ell_1,\ell_2):\,j\geq0\,,\ell_1,\ell_2=\pm2^j\,,k\in\Z^4\}$ boundary shearlets.
\end{itemize}
For $f \in L^2(\R^4)$, the {\it cylindrical shearlet transform} $\mathcal{S}$ is the mapping 
$$ f \mapsto \mathcal{S}(f)=\ip{f}{\tilde \psi_\mu}, \quad \mu \in M.$$

We remark that, by a direct computation, we can write the shearlet functions in~\eqref{sh.one} as 
\begin{equation}\label{shearTime}
\psi_{j,\ell,k}^{(d)}(x)=|\det A_{(d)}|^{j/2}\psi_{j,\ell}^{(d)}\left(B_{(d)}^{[\ell]}A_{(d)}^j x + k \right)
\end{equation}
where $\hat{\psi}_{j,\ell}^{(d)}(\xi)=W\left(2^{-2j}\xi B_{(d)}^{[\ell]}A_{(d)}^j\right)\!V_{(d)}(\xi)$ depends mildly on $j \ge 0$ and $\ell =(\ell_1,\ell_2) \in \Z^2$, where 
$|\ell_1|, |\ell_2| \le 2^j$.
Using the support and regularity of $W$ and $V$, one can show that, for any $\nu \in(\N\cup\{0\})^4$ and any $N>0$, there is $C_{\nu,N}>0$ independent of $j,\ell,d$ such that 
\begin{equation}\label{boundDeriv}
    \partial_{x}^{\nu}\psi_{j,\ell}^{(d)}(x)\leq C_{\nu,N}(1+|x|^2)^{-N}.
\end{equation}
The proof of this estimate is presented in Appendix~\ref{proof.ineq}.

\section{Sparse cylindrical shearlets approximations}   \label{s.main}

We start by defining the class of 4-dimensional cylindrical cartoon-like functions associated with our data model. This definition extends a similar definition in the 3-dimensional setting that was introduced by some of the authors~\cite{LPGE00} as a modification of the better known class of {\it cartoon-like functions}, originally proposed by Donoho \cite{Don01} to provide a simplified model of natural images. 

\subsection{Cylindrical cartoon-like functions}
\label{ssec:CylCartoonFunct}

For a fixed constant $A>0$, let $\M(A)$ be a class of indicator functions of sets $B\subset[0,1]^3$ with $C^2$-regular 2-manifold boundary $\partial B=\bigcup_\alpha \Sigma_\alpha$, where $\alpha$ ranges over a finite index set, and for each $\alpha$, the surface $\Sigma_\alpha$ has a parametrization $\Sigma_\alpha=\{(v,E_\alpha(v))\,:\,v\in V_\alpha\subset\R^2\}$, where  $E_\alpha$ is a $C^2$-regular function with values on the open set $V_\alpha \subset \R^2$, such that $\|E_\alpha\|_{C^2(V_\alpha)}\leq A$. Denoting with $C^2([0,1]^3)$ the collection of twice differentiable functions supported  inside $[0,1]^3$, we define the class of {\sl 4-dimensional cylindrical cartoon-like functions} $\E(A)$ as the set 
$$\{f=h_0g_0+h_1\X_B g_1\,:\, \X_B\in\M(A), \, h_0,h_1\in C^2\left([0,1]^3\right),\,g_0,g_1\in C^2([-1,1])\}$$
where 
\begin{equation}\label{cartLike}
    f(x_1,x_2,x_3,x_4)=h_0(x_1,x_2,x_3)g_0(x_4)+h_1(x_1,x_2,x_3)\X_B(x_1,x_2,x_3)g_1(x_4)
\end{equation}
and $\|f\|_{C^2}=\sum_{|\alpha|\leq 2}\|D^\alpha f\|_\infty\leq1$.

\subsection{Approximation theorems}

Let $\{\tilde{\psi}_\mu\}_{\mu\in M}$ be the the Parseval frame of cylindrical shearlets given by~\eqref{shearSystem}. The cylindrical shearlet coefficients of $f\in L^2(\R^4)$ are the elements of the sequence $\{s_\mu(f)=\ip{f}{\tilde{\psi}_\mu}:\, \mu\in M\}$. We denote by $|s_\mu(f)|_{(N)}$ the $N$-th largest entry in modulus of this sequence. 

We can now state our main theoretical result, whose proof is presented in Sec.~\ref{th31proof}. 

\begin{theorem}\label{mainTheo}
Let $f\in\E(A)$ and $\{s_\mu(f)=\ip{f}{\tilde{\psi}_\mu}:\, \mu\in M\}$ be the sequence of corresponding cylindrical shearlet coefficients. Then, for any $N\in\N$, there is a constant $C$ independent of $\mu$ and $N$ such that \begin{equation}\label{mainTheRes}
    \sup_{f\in\E(A)}|s_\mu(f)|_{(N)}\leq C N^{-1}(\log N).
\end{equation}
\end{theorem}

Let $f_N^S$ be the $N$-th term approximate of $f\in\E(A)$ obtained from the $N$-th largest coefficients of its cylindrical shearlet expansion, namely $f_N^S=\sum_{\mu\in I_N}\ip{f}{\tilde{\psi}_\mu}\tilde{\psi}_\mu$ where $I_N\subset M$
 is the set of indices corresponding to the $N$-th largest entries of the sequence $\{|\ip{f}{\tilde{\psi}_{\mu}}|^2:\,\mu\in M\}$. The approximation error satisfies the estimate:
 $$\|f-f_N^S\|_{L^2}^2\leq\sum_{m>N}|s_\mu(f)|_{(m)}^2.$$

 Thus,  Theorem~\ref{mainTheo} implies the following result directly.
 \begin{theorem} \label{mainTheo2}
 Let $f\in\E(A)$ and $f_N^S$ be the $N$-th term approximation defined above. Then, for $N\in\N$, there is a constant $C$ independent of $N$ and $\mu$ such that $$\|f-f_N^S\|_{L^2}^2\leq C N^{-1}(\log N)^2.$$
 \end{theorem}

{\bf Remark.} The decay estimate above is the same as the one found for 3-dimensional shearlets \cite{GL_3D} which is the optimal rate in the class of 3-dimensional cartoon-like functions \cite{GL_3D,kutyniok2012} and is faster than the optimal decay rate $O(N^{-2/3})$ valid for the class of 4-dimensional cartoon-like functions (cf.\cite{kutyniok2012}). We conjecture that $O(N^{-1})$ is indeed the optimal decay rate in the class of 4-dimensional cylindrical cartoon-like images. In particular, the decay rate of 4-d cylindrical shearlets is significantly faster than conventional 4-d wavelets whose decay rate is $O(N^{-1/3})$.

\subsection{Arguments and constructions}   \label{th31proof}
The general structure of the proof of Theorem~\ref{mainTheo} is similar to the structure of~\cite{GL_3D}. However, to deal with the geometry of 4d cylindrical shearlets, we need to introduce new technical constructions and modify some critical steps of the original arguments, especially in the proofs of Theorems~\ref{theo0} and \ref{theo1}  below.  

To measure the sparsity of shearlet coefficients, we introduce the weak-$\ell^p$ quasi-norm $\| \cdot \|_{w\ell^p}$ which,  for a sequence $s=(s_\mu)_{\mu\in M}$, is defined as $$\|s\|_{w\ell^p}=\sup_{N>0}N^{1/p}|s_\mu|_{(N)}$$
where $|s_\mu|_{(N)}$ is the $N-th$ largest entry in the sequence $s$. In~\cite{stein2016}, this norm is shown to be equivalent to $$\|s\|_{w\ell^p}=\left(\sup_{\epsilon>0} \#\{\mu:\, |s_\mu|>\epsilon\}\epsilon^p\right)^{1/p}.$$

We intend to analyze the decay properties of the cylindrical shearlet coefficients $\{\ip{f}{\tilde{\psi}_\mu}:\,\mu\in M\}$, where $f$ is chosen according to our cylindrical cartoon-like model~\eqref{cartLike}, that is,  $$f(x_1,x_2,x_3,x_4)=h(x_1,x_2,x_3)\X_B(x_1,x_2,x_3)g(x_4),$$
where $\X_B\in\M(A)$, $h\in C^2([0,1]^3)$, $g\in C^2([0,1])$.

 We recall that \textit{discontinuities only occur in the $x_1 x_2 x_3$ space}. Thus, to carry out our analysis, we smoothly localize the function $f$ near dyadic squares in the $x_1x_2x_3$ space as follows. For a scale parameter $j\geq0$ fixed, we let $$\Q_j=\left\{Q=\left[\f{k_1}{2^j},\f{k_1+1}{2^j}\right]\times\left[\f{k_2}{2^j},\f{k_2+1}{2^j}\right]\times\left[\f{k_3}{2^j},\f{k_3+1}{2^j}\right]:\, k_1,k_2,k_3\in\Z\right\}$$
be the collection of dyadic cubes. For a non-negative $C^\infty$ function with support $w$ in $[-1,1]^3$ we define a smooth partition of unity 
$$\sum_{Q\in\Q_j}w_Q(x)=1, \quad x\in\R^3,$$
where for each dyadic cube $Q\in\Q_j$, $w_Q(x)=w(2^{j}x-k)$ and $k\in\Z^3$. We will examine the cylindrical shearlet coefficients of $f_Q:=f w_Q$, i.e., $\{\ip{f_Q}{\tilde{\psi}_\mu}:\, \mu\in M_j\}$, where $M_j=\{\mu\in M: j \text{ is fixed}\}$. As we show below, these coefficients exhibit a different decay behaviour depending on whether the surface $\partial B$ intersects the support of $w_Q$ or not. Let $\Q_j=\Q_j^0\cup\Q_j^1$ be the disjoint union of $\Q_j^0=\{Q\in\Q_j:\,\partial B\cap \supp(w_Q)\neq \emptyset\}$ and $\Q_j^1=\{Q\in\Q_j:\,\partial B\cap \supp(w_Q)= \emptyset\}$. Notice each $Q$ has side-length $3\cdot 2^{-j}$, then $|\Q_j^0|\lesssim2^{2j}$. Similarly, since $\supp(f)\subset[0,1]^4$, then $\Q_j^1\lesssim 2^{3j}$.
With this notation, we now state two theorems that will be used to prove Theorem~\ref{mainTheo}. Note that while the decay rate in Theorem~\ref{theo0} is the same as the one found for 3-dimensional shearlets in \cite{GL_3D}, the decay rate in Theorem~\ref{theo1} is different.

\begin{theorem}\label{theo0}
 Let $f\in\E(A)$. For $Q\in\Q_j^0$, with $j\geq 0$ fixed, the cylindrical shearlet coefficients $\{\ip{f_Q}{\tilde{\psi}_\mu}:\, \mu\in M_j\}$ satisfy $$\|\ip{f_Q}{\tilde{\psi}_\mu}\|_{w\ell^1}\leq C 2^{-2j},$$
 where $C$ is a constant independent of $Q$ and $j$.
\end{theorem}

\begin{theorem}\label{theo1}
 Let $f\in\E(A)$. For $Q\in\Q_j^1$, with $j\geq 0$ fixed, the cylindrical shearlet coefficients $\{\ip{f_Q}{\tilde{\psi}_\mu}:\, \mu\in M_j\}$ satisfy $$\|\ip{f_Q}{\tilde{\psi}_\mu}\|_{w\ell^1}\leq C 2^{-7j/2},$$
 where $C$ is a constant independent of $Q$ and $j$.
\end{theorem}

We show next how to apply Theorems~\ref{theo0} and \ref{theo1} to prove Theorem~\ref{mainTheo}
using an argument  similar to~\cite{GL_3D}; we postpone their rather technical proofs to the Appendix~\ref{s.proof}.

We have the following corollary. 
\begin{corollary}\label{mainCor}
Let $f\in\E(A)$ and for $j\geq 0$, consider the sequence of cylindrical shearlets coefficients $s_j(f)=\{\ip{f}{\tilde{\psi}_\mu}:\, \mu\in M_j\}$. Then there is a constant $C$ independent of $j$ such that $$\|s_j(f)\|_{w\ell^1}\leq C.$$
\end{corollary}
\textbf{Proof.} Using Theorems~\ref{theo0}, \ref{theo1} and the triangle inequality for weak $\ell^1$ spaces, we have
\begin{eqnarray*}
\|s_j(f)\|_{w\ell^1}&\leq&\sum_{Q\in\Q_j}\|\ip{f_Q}{\tilde{\psi}_\mu}\|_{w\ell^1}\\
&\leq&\sum_{Q\in\Q_j^0}\|\ip{f_Q}{\tilde{\psi}_\mu}\|_{w\ell^1}+\sum_{Q\in\Q_j^1}\|\ip{f_Q}{\tilde{\psi}_\mu}\|_{w\ell^1}\\
&\leq& C|\Q_j^0|2^{-2j}+C|\Q_j^1|2^{-7j/2}\\
&\leq& C(2^{2j}2^{-2j}+2^{3j}2^{-7j/2})\leq C. 
\end{eqnarray*}
In the last step, we have used the observations that $|\Q_j^0|\leq C 2^{2j}$ and $|\Q_j^1|\leq C 2^{3j}$. $\qed$

We next prove Theorem~\ref{mainTheo}.

\textbf{Proof of Theorem~\ref{mainTheo}.} By Corollary~\ref{mainCor}, we have that
\begin{equation}\label{numElements}
    R(j,\epsilon)=\#\{\mu\in M_j:\,|\ip{f}{\tilde{\psi}_\mu}|>\epsilon\}\leq C\epsilon^{-1}
\end{equation}
For an interior shearlet $\psi_{j,\ell,k}^{(d)}$, given by~\eqref{sh.one.f}, a direct computation using~\eqref{shearTime} and~\eqref{boundDeriv} gives that
\begin{eqnarray}
|\ip{f}{\psi_{j,\ell,k}^{(d)}}|&=&\left|\int_{\R^4}f(x) |\det A_{(d)}|^{j/2}\psi_{j,\ell}^{(d)}\left(B_{(d)}^{[\ell]}A_{(d)}^j x + k \right) dx\right|\nonumber\\
&\leq&2^{3j}\|f\|_\infty\int_{\R^4}\left|\psi_{j,\ell}^{(d)}\left(B_{(d)}^{[\ell]}A_{(d)}^j x + k \right)\right|dx\nonumber\\
&\leq&2^{-3j}\|f\|_\infty \int_{\R^4}|\psi_{j,\ell}^{(d)}(y)|dy\nonumber\\
&\leq&C2^{-3j}. \label{coefBound}
\end{eqnarray}
A very similar computation on the boundary shearlets gives the same estimate.
So, for a given $\epsilon>0$, there is $j_\epsilon>0$ such that $|\ip{f}{\tilde{\psi}_{\mu(j)}}|<\epsilon$ for each $j\geq j_\epsilon$. Therefore, from~\eqref{coefBound}, we have $R(j,\epsilon)=0$ for $j>\f{1}{3}\log_2(\epsilon^{-1})+\log_2(C)>\f{1}{3}\log_2(\epsilon^{-1})$. So, using~\eqref{numElements}, we have
\begin{equation*}\label{numBound}
\#\{\mu\in M:\,|\ip{f}{\tilde{\psi}_\mu}|>\epsilon\}\leq\sum_{j\geq0}R(j,\epsilon)=\sum_{j=0}^{\f{1}{3}\log_2(\epsilon^{-1})}R(j,\epsilon)\leq C\epsilon^{-1}\log_2(\epsilon^{-1}).
\end{equation*}
Next, let $n=n(\epsilon) = \#\{\mu\in M:\,|\ip{f}{\tilde{\psi}_\mu}|>\epsilon\}$. Notice $\epsilon^{-1}\lesssim n$. Therefore, from~\eqref{numBound}, we have $\epsilon\leq C n^{-1} \log_2(\epsilon^{-1})\leq C n^{-1} \log_2(n)$. So, if $|s(f)|_{(N)}$ is the $N$-th largest coefficient, then $|s(f)|_{(N)}\leq C N^{-1} \log_2(N)$ and inequality \eqref{mainTheRes} follows. $\qed$


\section{Numerical implementation of 4d cylindrical shearlet}
\label{s.numImplementation}

This section covers the practical implementation of the 4d cylindrical shearlet transform, including its inverse and adjoint transforms, which we apply in Sec.~\ref{s.numDynamicCT} to illustrate the potentiality of cylindrical shearlets in numerical applications. Numerical codes for the Matlab framework, with all the necessary documentation, are available in Github \cite{heikkila20214dCylSh}.

\begin{figure}[bth!]
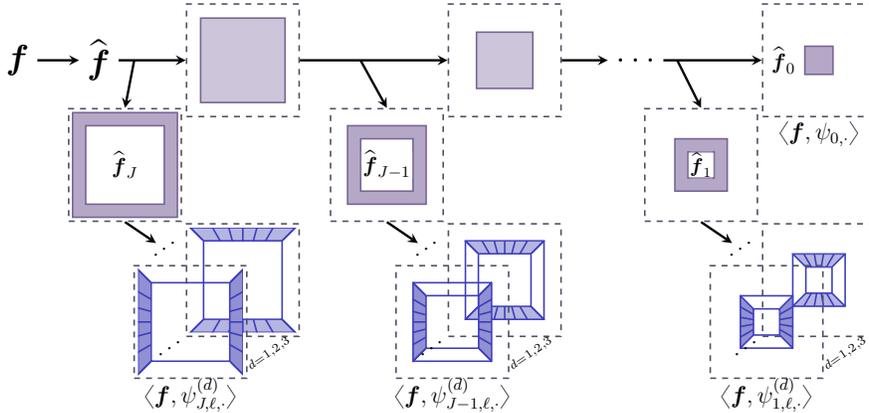

    \centering
    \addPDFfig{1}
    \caption{Illustration of how the cylindrical shearlet transform decomposes an input signal.  The signal is split into $J+1$ subbands in the Fourier domain using the Laplacian pyramid scheme (step 1). Next, cylindrical shearlet coefficients are computed  for all scales, orientations and hyperpyramids (step 3) using directional filters associated with the three hyperpyramids and different shearing parameters $\ell = (\ell_1, \ell_2)$ (created in step 2).}
    \label{fig:decompScheme}
\end{figure}

Key ideas of the implementation generalize those of the 3d cylindrical shearlet transform \cite{LPGE00} and are illustrated by the decomposition scheme in Fig.~\ref{fig:decompScheme}. The directional filters, in particular, are derived from the 3d discrete shearlets \cite{NL2012} since the transform only captures directional structures along the first three dimensions. We point out though that the current implementation is meant as a proof of concept and possible solutions to the evident inefficiencies are left for future work. 

In the following, we denote discrete values at a specific multi-index by using square brackets $[\cdot]$. Operations like the inner product $\langle \cdot, \cdot \rangle$, the discrete Fourier transform $\F (x)[\xi] = \widehat{x}[\xi]$ and its inverse $\F^{-1} (x) = \widecheck{x}$ are defined as usual, unless otherwise specified.
Given a discretized 4d object $\fvec$ and $J,\ell,k$ fixed (notice that these indexes depend on the resolution and user inputs), the forward transform $\CylShTD$ can be computed using the following steps:
\begin{enumerate}
    \item \textbf{Subband decomposition.} Compute the 4d multilevel (up to scale $J$) subband decomposition $\P_J\fvec = \left(\fvec_j\right)_{j=0}^J$ using an adapted Laplacian pyramid scheme\footnote{This particular implementation of the Laplacian pyramid decomposition was originally introduced for the surfacelet transform \cite{lu2007} and here we follow~\cite{LPGE00} in applying the same idea to cylindrical shearlets.}. In the frequency domain, each subband corresponds to a windowing:
    \[
        \widehat{\fvec}_j[k] = \widehat{\fvec}[k] \Wvec_j[k].
    \]
    This operation corresponds to the first part of equation~\eqref{shear1} in the construction of the Parseval frame (see Sec.~\ref{ss.parsevalConstruction}).
    \item \textbf{Directional filtering.} Consistently with Definition~\ref{def:cylindrical} and similar to~\cite{NL2012}, we construct the directional filters $\Vvec_{j,\ell}^{(d)}$, for given scale $j$, direction $\ell = (\ell_1, \ell_2)$ and cylindrical hyperpyramid indexed by $d$, by first defining a window in the pseudo-spherical Fourier domain and then resampling on a Cartesian grid. Note that directional filters are (approximately) symmetric and, by construction, the decomposed signal does not need resampling. 
    \item \textbf{Transform coefficients.} Using the directional filters, the cylindrical shearlet coefficients are thus given by:
    \begin{equation} \label{eq.numCoef}
        \CylShTD (\fvec) := \langle \fvec, \psi^{(d)}_{j,\ell,k} \rangle = \fvec_j \cc \widecheck{\Vvec}_{j,\ell}^{(d)} [k] = \F^{-1}\left( \widehat{\fvec}_j [k]\Vvec_{j,\ell}^{(d)}[k] \right),
    \end{equation}
    where $k = (k_1, k_2, k_3, k_4) \in \Z^4$ and $\cc$ denotes the discrete convolution along the first three dimensions.
\end{enumerate}

In step 2, the initial construction of the $\Vvec_{j,\ell}^{(d)}$'s does not guarantee that the sizes match the size of $\widehat{\fvec}_j$. Hence, the necessary padding and many FFTs required slow down the algorithm considerably. In our experiments to perform step 3, we computed the cylindrical (3d) convolution as a pointwise multiplication in the Fourier domain to reduce repetitive operations.


From the cylindrical shearlet coefficients~\eqref{eq.numCoef}, the original 4d signal $\fvec$ is recovered using the inverse cylindrical shearlet transform $\CylShTD^{-1}$. By design of the directional filters, for each scale $j$, we have: 
\begin{equation*}
    \sum_{d=1}^3 \sum_{\ell} \Vvec_{j,\ell}^{(d)} = 1.
\end{equation*}
Hence, the inverse transform is very straightforward to compute since we have:
\begin{eqnarray*}
    \sum_{d=1}^3 \sum_{\ell} \langle \fvec, \psi^{(d)}_{j,\ell,\cdot} \rangle = \sum_{d=1}^3 \sum_{\ell} \F^{-1}\left( \widehat{\fvec}_j [k]\Vvec_{j,\ell}^{(d)}[k] \right) = \fvec_j,
\end{eqnarray*}    
which in turn yields
\begin{equation} \label{eq.inverseTransform}
    \fvec = 
\CylShTD^{-1}\left( \langle \fvec, \psi^{(d)}_{j,\ell,k} \rangle \right) := 
\P_J^{-1} \left(\sum_{d=1}^3 \sum_{\ell} \langle \fvec, \psi^{(d)}_{j,l,\cdot} \rangle \right),
\end{equation}
where $\P_J^{-1}$ is the inverse of the Laplacian pyramid-like scheme. This computation is very efficient as no convolutions or filters are needed.

Finally, we discuss how to implement the adjoint (or synthesis) operator $\CylShTD^*$ whose computation is more involved than the inverse. First we let $\uShc[k]$ be a vector on the cylindrical shearlet coefficient domain. Then a direct computation of the inner product gives:
\begin{eqnarray*}
    \big\langle \CylShTD(\fvec), \uShc \big\rangle &=& \sum_{j=0}^J \sum_{d=1}^3 \sum_{\ell} \big\langle \langle \fvec, \psi^{(d)}_{j,\ell,\cdot} \rangle , \uShc \big\rangle \\
    &=& \sum_{j=0}^J \sum_{d=1}^3 \sum_{\ell} \langle \fvec_j \cc \widecheck{\Vvec}_{j,\ell}^{(d)}, \uShc \rangle \\
    &=& \sum_{j=0}^J \Big\langle \fvec_j, \sum_{d=1}^3 \sum_{\ell} \uShc \cc \vec{\Lambda}_{j,\ell}^{(d)} \Big\rangle \\
    &=& \Big\langle \fvec, \P_J^* \left(\sum_{d=1}^3 \sum_{\ell} \uShc \cc \vec{\Lambda}_{j,\ell}^{(d)} \right) \Big\rangle \\
    &=& \Big\langle \fvec, \CylShTD^{-1} \left(\uShc \cc \vec{\Lambda}_{j,\ell}^{(d)} \right) \Big\rangle,
\end{eqnarray*}
where $\vec{\Lambda}_{j,\ell}^{(d)}[k_1, k_2, k_3, k_4] := \widecheck{\Vvec}_{j,\ell}^{(d)}[-k_1, -k_2, -k_3, k_4] \approx \widecheck{\Vvec}_{j,\ell}^{(d)}[k_1, k_2, k_3, k_4]$ by symmetry. The last step follows from equation~\eqref{eq.inverseTransform} and the observation that, unlike a traditional Laplacian pyramid decomposition (cf.~\cite{BA1987}), the implementation in \cite{lu2007} has the property that $\P_J^* = \P_J^{-1}$; this simplifies the final implementation step. However, computing the adjoint is slower than the inverse and the end result is slightly blurred due to the convolutions involved.

\section{An application to dynamic tomography}
\label{s.numDynamicCT}
In this section, we illustrate the numerical advantages of cylindrical shearlets vs.~conventional 4d wavelets when dealing with spatio-temporal data by considering a challenging inverse problem, namely the reconstruction of a volume over time associated with 4d (3d+time) dynamic CT.

CT is a classical inverse problem concerned with recovering the inner structure of an unknown object from external measurements of its X-ray attenuation intensity. This task is notoriously ill-posed, especially when measurements are sparse. One way to overcome ill-posedness and, thus, to guarantee a stable and unique solution, is to add \textit{regularization} to the problem~\cite{Engl1996}. During the last decade, several sparse regularization strategies 
were proposed in CT applications, based on the paradigm that, for any data class, there exists an appropriate sparsifying data representation, e.g., wavelets or shearlets.

Here we illustrate the application of a regularizer based on cylindrical shearlets to dynamic CT
by adapting to the 3d+time setting a regularized reconstruction method based on (conventional) shearlets proposed by some of the authors in~\cite{bubba2020sparse}. This reconstruction method was originally motivated by sparse imaging of phloem transport in plant stems and was shown to be extremely competitive as compared to other methods from the literature.

A main advantage of this approach is that, unlike many existing methods~\cite{Katsevich10,Roux04}, is not limited to 2d data, and, unlike methods relying on filtered back-projection (FBP), cf.~\cite{blanke2020inverse,Hahn14,Hahn16,HahnQuinto16}, does not require a dense angular sampling. Additionally, we do not need to assume periodicity on the movement as in~\cite{Ritman03}, nor constant total brightness as compared to optical flow~\cite{burger2017} nor multiple source–detector pairs, as in~\cite{Hakkarainen19,Niemi15}. 
We refer the interested reader to~\cite{Hauptmann2021} for a broader overview of image reconstruction in dynamic inverse problems.

\subsection{Mathematical model} \label{ss.mathematicalModel}

Modern cone-beam CT scanners reconstruct a 3d volume of the interior attenuation of the targeted object using 2d projection images collected from multiple angle views. If this measurement process is repeated over time, the object of interest can be understood as a 4d object. As observed above, given the sparse measurements and the violation of the static assumption that is often assumed in classical CT reconstruction schemes, stable recovery of a moving object from multiple sparse measurements over a time period requires regularization.

The novelty here is that, by applying 4d cylindrical shearlets for regularization, we do not only regularize over the 3d spatial volume but also across time frames within the same representation system. This property is expected to be a significant advantage with respect to separable representations due to the superior approximation properties of cylindrical shearlets that were discussed in Sec.~\ref{s.main}. This improved behavior is confirmed by our numerical results.

Formally, for each time step $t = 1,...,\tau$, let $\fvec_t[x_1,x_2,x_3] \in \R^{n}_+$, with $n = n_{x_1}n_{x_2}n_{x_3}$, be a vector representing the unknown 3d object, $\RadonD_t \in \R^{p \times n}$ a matrix modelling the tomographic cone-beam measurement process and $\m_t + \etab =: \m^{\etab}_t \in \R^p$ the data corrupted by measurement errors $\etab = \etab(t)$. To further simplify our notation we set:
\begin{equation*}
    \fvec = \left[ \begin{array}{c}
    \fvec_1 \\
    \vdots \\
    \fvec_{\tau} \end{array} \right], \ \RadonD = \left[ \begin{array}{ccc}
    \RadonD_1 & & \\
     & \ddots & \\
     & & \RadonD_{\tau} \end{array} \right], \ \m^{\etab} = \left[ \begin{array}{c}
    \m^{\etab}_1 \\
    \vdots \\
    \m^{\etab}_{\tau} \end{array} \right].
\end{equation*}
Then a regularized solution $\fvec \in \R_+^{n\tau}$ is obtained by minimizing the functional
\begin{equation} \label{eq:functional}
    J(\fvec) = \frac{1}{2} 
    \big\|\RadonD \fvec - \m^{\etab} \big\|_2^2 
    + \beta \|\CylShTD \fvec\|_1.
\end{equation}
Here, the regularization parameter $\beta > 0$ balances between the data mismatch term over the time steps and the $\ell^1$-sparsity of 4d cylindrical shearlet coefficients of the solution.

A robust minimization method is the Primal-Dual Fixed Point (PDFP) algorithm~\cite{chen2016primal}, which generalizes the well-known Iterative Soft-Tresholding Algorithm (ISTA) to include non-negativity constraints for the solution $\fvec$ and ensures convergence even when the sparsifying system does not form an orthonormal basis but a frame, which is the case with cylindrical shearlets as shown in Sec.~\ref{ss.parsevalConstruction}. By using PDFP, equation~\eqref{eq:functional} can be minimized by iterating the following steps:
\begin{equation} \label{eq:PDFP} 
\begin{cases}
 \y^{(i+1)} &= \proj{+} \big(\fvec^{(i)} - \rho (\RadonD^T \RadonD \fvec^{(i)} - \RadonD^T \m^{\etab}) - \lambda \CylShTD^* \rvec^{(i)} \big),  \\[0.25em]
\rvec^{(i+1)} &= \big(\mathbb{I} - S_{\beta \frac{\rho}{\lambda}}\big) \big(\CylShTD \y^{(i+1)} + \rvec^{(i)} \big), \\[0.25em]
\fvec^{(i+1)} &=\proj{+} \big(\fvec^{(i)} - \rho (\RadonD^T \RadonD \fvec^{(i)} - \RadonD^T \m^{\etab}) - \lambda \CylShTD^* \rvec^{(i+1)} \big) 
\end{cases}
\end{equation}
where $S_{\beta \frac{\rho}{\lambda}}$ denotes the soft-thresholding operator and $\proj{+}$ is the projection onto the non-negative orthant.
The parameters $\rho$ and $\lambda$ are bounded by properties of the functional $J$, which set a clear range for their values, while the optimal choice of $\beta$ is a notoriously difficult task.

Here, we adopt an automated tuning of $\beta$ based on the given a priori sparsity level of the cylindrical shearlet coefficients. This method was originally introduced in~\cite{purisha2017controlled} using Haar wavelet regularization in traditional 2d tomography regularization. 

In a recent work, some of the authors modified this method for the 2d+time dynamic tomography setting using classical shearlets~\cite{bubba2020sparse} and 3d+time (complex) wavelets~\cite{Bubba2021}, where they also provided further justification for this model. The detailed steps of this method are found in Algorithm 1 in~\cite{bubba2020sparse}, where the necessary modifications from 3d to 4d apply.

\subsection{Simulated test data} \label{ss.data}
Our regularized reconstruction approach is applied to a simulated 4d tomography dataset consisting of repeated measurements of a custom ellipsoid phantom created using the 3d phantom from~\cite{jorgensen2010tomobox}. The intensity values of the two larger ellipsoids change linearly in the interval $[0,1]$ while the intensities of the multiple smaller ellipsoids follow a sinusoidal pattern with offset phases. The codes for generating the data are included in the Github repository~\cite{heikkila20214dCylSh}.  We remark that the phantom is  consistent with our model of cylindrical cartoon-like functions as the spatial boundaries remain fixed.

\begin{figure}[bth!]
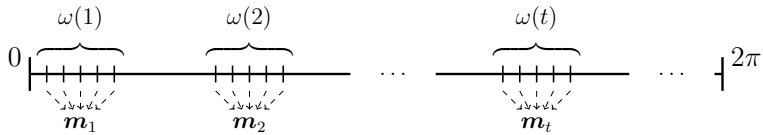

    \centering
    \addPDFfig{2}
    \caption{Illustration of how multiple values of $\omega$ are used to simulate measurements at each time frame $t$.}
    \label{fig:dynamicSimulation}
\end{figure}

The spatial dimensions of the volume are $256\times256\times64$ voxels\footnote{While seemingly small, storing this 4d object in single precision already requires roughly 260 MB.} and we simulated in total $16$ sparse angle cone-beam sinograms (i.e., corresponding to 16 different time frames).
We tested varying number of evenly spaced projection angles: in Sec.~\ref{ss.results}, we report results with 24, 30, 60 and 90 equispaced angles. The matrices $\RadonD_t$ (and therefore $\RadonD$) simulating the geometry of a cone-beam CT are generated using the HelTomo Toolbox~\cite{meaney2021heltomo}, build upon the ASTRA Toolbox~\cite{van2016fast}. 
All changes depend on the parameter $\omega \in [0,2\pi]$ and to better approximate continuous motion and realistic measurement conditions we sample it as follows. First, given the number of desired measurements (here $\tau = 16$), we divide the whole period $[0,2\pi]$ into $2\tau - 1$ subintervals and discard every second one leaving a total of $\tau$ disjoint subintervals. Then, for each $t$, we further sample the corresponding subinterval to obtain 15 values of $\omega(t) = (\omega_1^{(t)},\dots, \omega_{15}^{(t)})$ which are used to simulate a sinogram $\m_t$ in 15 stages. The middle value from each subinterval (i.e., $\omega_8^{(t)}$ for each $t$) is taken as the ground truth $\fvec_t$ to benchmark the reconstructions. The sampling procedure is illustrated in Fig.~\ref{fig:dynamicSimulation}.

\begin{figure}[tbh!]
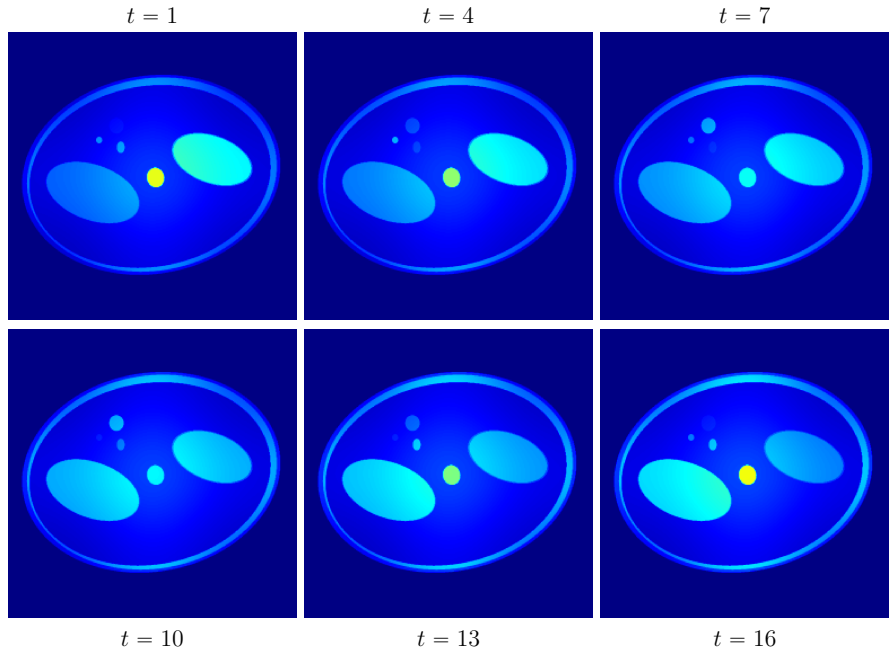

	\centering
    \addPDFfig{3}
	\caption{Interior ($xy$-plane, $z=28$) slices of the phantom as it evolves over time steps $t = 1,4,7,10,13,16$.}
	\label{fig:CPphantom}
\end{figure}

The jump in values of $\omega$ between measurements simulate a pause between consecutive measurement cycles and allows for more noticeable changes between time frames, as illustrated in Fig.~\ref{fig:CPphantom}, which contains several horizontal ($xy$-plane) slices of the phantom at different time frames $t$ covering the full range of changes. In comparison, changes allowed within each $t$ are less severe but not negligible as in practice the measurement device cannot measure all projections simultaneously.

Finally, to avoid inverse crime, each individual projection image is generated at twice the desired resolution, down-sampled and then corrupted by white Gaussian noise (0 mean and 5\% variance).

\subsection{Results} \label{ss.results}
We report here the numerical results of our reconstruction algorithm for dynamic CT. Results comprise reconstructions from a varied number of sparse projection angles. For simplicity the same evenly spaced angular sampling is used for every time frame.

For comparison purposes, we implemented also a regularized reconstruction algorithm based on 
the 4d discrete wavelet transform (DWT). The regularized model with 4d DWT is obtained by replacing $\CylShTD$ with a DWT in equation~\eqref{eq:PDFP} and changing the values of $\lambda$ and $\beta$ accordingly. The 4d DWT is implemented by extending the 3d DWT from Matlab's Wavelet Toolbox and it is available on GitHub~\cite{heikkila2021wavedec4}. The wavelet decomposition, based on Daubechies 2 filters, is performed using 4 scales. The cylindrical shearlet transform uses 3 scales with the number of directions in each pyramid being 36, 16 and 4 as the scale gets coarser. We found this setting to offer a good balance between quality and efficiency.

Due to the high memory requirements of the 4d cylindrical shearlet transform, the computations were carried out on the Turso cluster at the University of Helsinki, using $16$ CPU cores each equipped with $16$GB memory. The 4d wavelet computations were performed on the same cluster but required just $8$GB of memory in total.

\begin{figure}[tbh!]
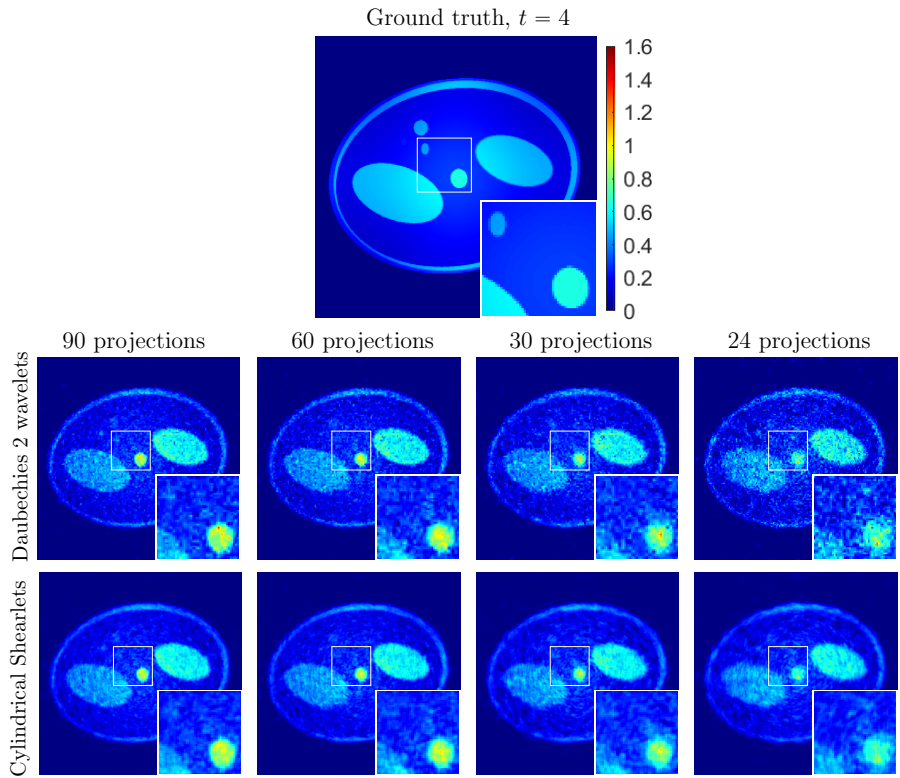

	\centering
    \addPDFfig{4}
    \caption{Interior slices at time frame $t = 4$ of the ground truth (top row), wavelet-regularized solution (center row) and cylindrical shearlet regularized solution (bottom row). On different columns there are reconstructions from 90, 60, 30 and 24 projections.}
\label{fig:t4recn}
\end{figure}

\begin{figure}[tbh!]
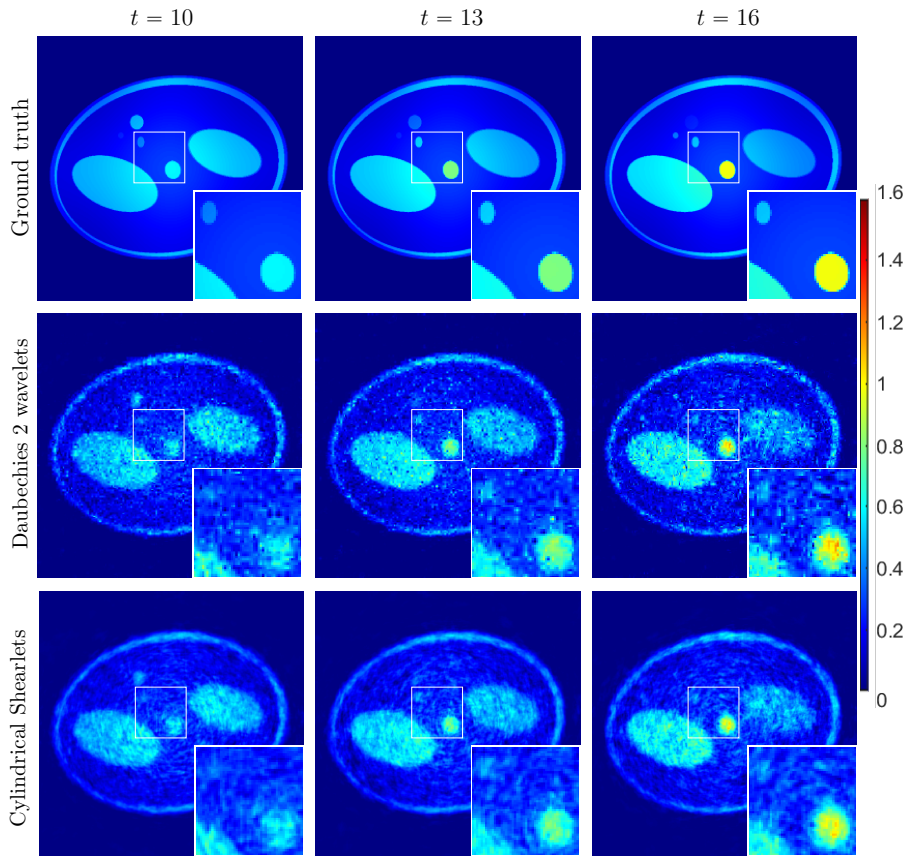

	\centering
    \addPDFfig{5}
    \caption{Interior slices of the ground truth (top row), wavelet-regularized solution (center row) and cylindrical shearlet regularized solution (bottom row) evolving over time frames $t = 10$ (left), $t=13$ (middle) and $t=16$ (right). All reconstructions are from 30 projections.}
\label{fig:30anglesrecn}
\end{figure}

\begin{figure}[tbh!]
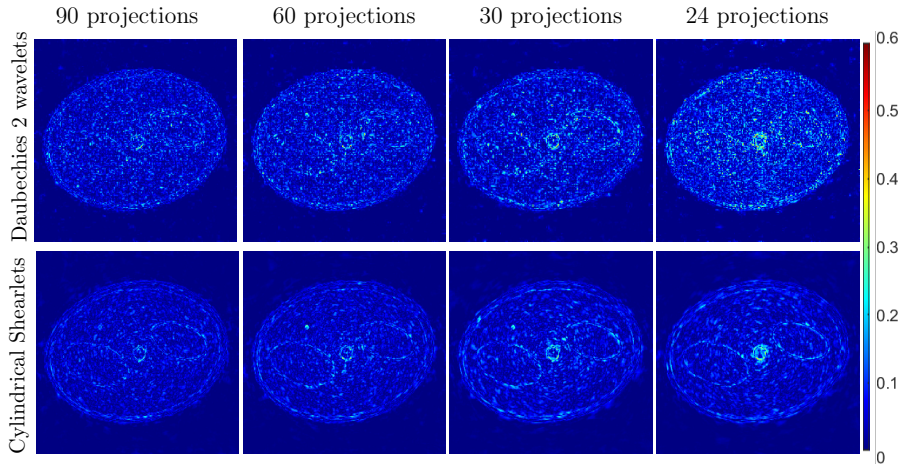

    \centering
    \addPDFfig{6}
    \caption{Absolute difference images $| \fvec_{\text{recn}} - \fvec_{\text{true}} |$ at time frame $t = 4$ of the wavelet regularized solution (top row) and cylindrical shearlet regularized solution (bottom row). On different columns there are reconstructions from 90, 60, 30 and 24 projections. Note that the color scaling is different in this figure.}
    \label{fig:t4dif}
\end{figure}

In Fig.~\ref{fig:t4recn} and Fig.~\ref{fig:30anglesrecn}, we display selected interior slices of our phantom reconstructed using the proposed algorithm and either 4d wavelets or 4d cylindrical shearlets for regularization. Fig.~\ref{fig:t4recn} displays a varying number of projections and the ground truth at a fixed time frame ($t = 4$). Fig.~\ref{fig:30anglesrecn} displays multiple time frames including the ground truth with the number of projections fixed to 30. To highlight the differences between the two regularization approaches, we display in Fig.~\ref{fig:t4dif} the absolute difference between the reconstructions and the true objective. In this figure, the time frame is again fixed at $t=4$. 

In all figures, we only display an $xy$-plane of the reconstruction at height $z=28$. Additionally, we include as an insert a zoomed-in sub-region containing smaller key details to highlight the reconstruction quality at the discontinuities. Reconstructions along other cross-sections of the solid exhibit similar properties. Colors are scaled uniformly so that all images in Figures \ref{fig:t4recn} and \ref{fig:30anglesrecn} (and \ref{fig:CPphantom}) are comparable. The color values in Fig.~\ref{fig:t4dif} are only comparable within that figure.

We report numerical error metrics in Table~\ref{tab:errors}. Specifically, we computed the Peak-Signal-to-Noise-Ratio (PSNR) comparing the whole 4d reconstruction to the known ground truth. In addition, we used the recently introduced Haar-wavelet Perceptual Similarity Index (HPSI)~\cite{reisenhofer2018}, originally proposed for images (i.e., 2d data) and here adapted to handle our higher dimensional data: namely, for each value of $t$, we compare the central slice to the central slice of the ground truth and report the mean value across all time frames. Finally, we compute the 3d Structural SIMilarity index (SSIM) \cite{wang2004image}
which we also average across all time frames. This last error metric is arguably the most faithful among those in Table~\ref{tab:errors} as it considers the whole 4d data (unlike HPSI) and, to an extent, the geometry (unlike PSNR).

\begin{table}[t]
    \centering
    \caption{Numerical error estimates of the different reconstructions.}
    \label{tab:errors}   
    \begin{tabular}{l|c c c c}
        & \multirow{2}{*}{Projections} & \multirow{2}{*}{PSNR} & Mean & Mean \\ 
        &  &  & HPSI & SSIM \\\hline
        \multirow{4}{2.5cm}{Daubechies~2-wavelets} & 90 & 29.0 & 0.500 & 0.815 \\
        & 60 & 28.3 & 0.467 & 0.798 \\ 
        & 30 & 27.6 & 0.429 & 0.782 \\
        & 24 & 25.3 & 0.381 & 0.716 \\[0.5em]
        \multirow{4}{2.5cm}{Cylindrical shearlets} & 90 & 31.1 & 0.576 & 0.841 \\
        & 60 & 30.6 & 0.547 & 0.826 \\
        & 30 & 29.8 & 0.500 & 0.801 \\
        & 24 & 29.4 & 0.487 & 0.789
    \end{tabular}
\end{table}

Results show that overall regularization based on cylindrical shearlet  yields a better reconstruction performance both in terms of visual quality and quantitative performance metrics. Visual differences in the reconstructions
are more pronounced for sparse projections (see Fig.~\ref{fig:30anglesrecn}, and the two rightmost columns of Fig.~\ref{fig:t4recn}) than for denser ones (Fig.~\ref{fig:t4recn}, left columns). In particular, cylindrical shearlet reconstructions are consistently better at suppressing noise without producing excessive blur near edges while wavelet regularized solutions are notably noisier. Indeed, wavelet-based reconstructions suffer from salt-and-pepper-like noise where the attenuation level is strongly under- or overestimated in some points (see for example center row of Fig.~\ref{fig:t4recn} or \ref{fig:t4dif}). The outer shell of the phantom is relatively well reconstructed by both regularized approaches.

Quantitative error measures in Table~\ref{tab:errors} show that cylindrical shearlets are better performing under all the error metrics. Even with just 30 projections, the cylindrical shearlet reconstruction outperforms the best wavelet reconstruction (90 projections) based on most of these metrics. We explain this behavior with the superior approximation properties of cylindrical shearlets, as shown above in Theorem~\ref{mainTheo}, that are manifested by their improved noise suppression as compared to wavelets. We also observe that cylindrical shearlets-based reconstructions are highly consistent across the number of projections; that is, reducing the number of projections does not result in a significant worsening of reconstruction  (the PSNR drops less than 2dB as the number of projection goes from 90 to 24) and visual quality. By contrast, wavelet-based reconstructions degrade rapidly as the number of projects decreases (the PSNR drops by over 2 dBs as the number of projections  go from 24 to 30) and visual inspection shows the appearance of salt-and-pepper-like artefacts.

Finally, we remark that while none of our error metrics fully accounts for the geometry of 4d spatio-temporal data, SSIM is able to assess in some form the overall reconstruction quality of 3d moving
volumes. Indeed, SSIM is applied in the literature to measure video approximations (2d + time setting) \cite{wang2004video}. By contrast, PSNR only accounts for  pointwise values with no geometric considerations and HPSI can be computed from 2d slices only (no extension to higher dimension is currently available). In the (medical) literature, the quality of 4d CT is often assessed using different correlation-type metrics (cf. \cite{castillo2014, noid2017}) which do not consider the geometry either. The topic of image and video quality assessment is a vast and active research area (cf. \cite{bampis2018spatiotemporal, chow2016review}); unfortunately, no fully satisfactory quality metric for 3d + time data is currently available.

\section{Discussion and conclusion} \label{s.conclusion}

We have introduced a new construction of multiscale representations on $L^2(\R^4)$ that is especially designed for the efficient approximations of spatio-temporal data. Our theoretical analysis shows that this method provides highly sparse representations in the class of 4d-cartoon-like images, outperforming conventional multiscale representations. We have also illustrated the practical advantages of the new representation on a challenging computational problem of regularized reconstruction in dynamic tomography from a small number of projections. Our numerical results show that our regularized reconstruction based on cylindrical shearlets outperforms a similar algorithm based on wavelets both in terms of visual quality and quantitative performance metric when projections are sparse. While our result was demonstrated using simulated data, we expect that a comparable performance advantage will hold using experimental data and will be investigated in a future work by extending our study of phloem transport in plant stems~\cite{bubba2020sparse}. In fact, the results of our study suggest a number of theoretical extensions and numerical refinements to further exploit the potential of cylindrical shearlets in numerical applications.

As mentioned above, the model of cylindrical cartoon-like functions adopted in this paper is a rather crude simplification of temporal sequences of 3-dimensional images found in applications, as it does not allow discontinuities with respect to the temporal variable. While the phantom we used in our simulations was designed to fit this model, realistic applications of dynamic CT typically involve boundaries in the spatial domain that change in time so that our image model would need to be modified.

\begin{figure}
    \centering
    \includegraphics[scale=0.5]{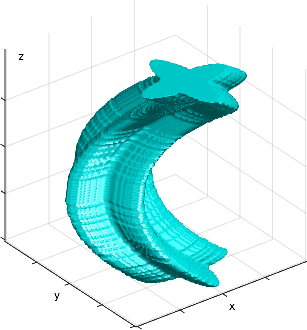}
    \caption{Modified 3-dimensional cartoon-like image. The discontinuous boundary curve in the $xy$ plane is rigidly displaced as a function of the $z$ coordinate.}
    \label{fig:model2}
\end{figure}

We are confident that the proofs presented in this work can be extended with a relatively simple argument to include generalizations of the cylindrical cartoon-like model such as the situation of a moving solid object, e.g., a moving ball. In this case, the boundary of the object is (smoothly)  displaced from a time-frame to the next one, without changes in the discontinuous boundary other than its location being rigidly translated (see Fig.~\ref{fig:model2} for an illustration in $\R^3$). To model such functions on $\R^4$, we may consider a modified cylindrical image model on $L^2(\R^4)$ where $f(x_1, x_2, x_3,x_4)= h((x_1,x_2,x_3)-t(x_4))g(x_4)$ where $h$ is a compactly supported $C^2$ function away from $C^2$ boundaries, $g \in C^2([-1,1]$ and $t$ is a smooth translation function which depends on $x_4$ only. The Fourier transform of $f$ is of the form $\widehat{f}(\xi_1,\xi_2,\xi_3,\xi_4)=\widehat{h}(\xi_1,\xi_2,\xi_3) \,  G(\xi_1,\xi_2,\xi_3,\xi_4)$ where $G$ is smooth and bounded. From this observation, it follows that one can adapt essentially the same arguments presented above to derive a result similar to Theorems~\ref{mainTheo} and \ref{mainTheo2}. A rigorous discussion of this extension of our proof would require more technical details that are beyond the scope of this paper and, for reasons of space, are left to a future work.

A downside of the current implementation of 4d cylindrical shearlets is the higher computational cost compared to 4d wavelets, due to the increased complexity of the transform. For the data size we considered, the computational burden is noticeable: the algorithm walltime for cylindrical shearlets is approximately 21 hours in total or 650 seconds per iteration, while for the wavelets total computing time is 1.5 hours or 22 seconds per iteration. A single computation of the forward, inverse or adjoint transform is manageable even on a regular desktop computer but iterative schemes usually require applying this computation hundreds of times, making the whole procedure very time consuming. On the other hand, the different number of X-ray projections (24, 30, 60 or 90) have a relatively minor impact on the overall computing time which is dominated by the cost of applying the wavelet or shearlet forward and adjoint transforms.
Nonetheless, we remark that our implementation of 4d cylindrical shearlets is
presented here as a proof of concept without a systematic effort to optimize the computational cost which would be beyond the scope of this paper. However, the parallel nature of the transform indicates a potentially significant speedup, for example, by utilizing GPU computing\cite{gibert2014discrete, andrade2020tfshearlab}.

Finally, we recall that deep learning strategies have gained increasing popularity in inverse problems including CT, where they have been applied very successfully often in combination with model-based principles such as sparsity models (e.g., \cite{bubba2020sparse}). We expect that the ideas presented in this paper have also the potential of being successfully integrated into a deep learning strategy leading to a new generation of reconstruction algorithms for dynamic CT integrating learning- and model-based principles.

\section*{Acknowledgements}
All authors acknowledge the support of the IT for Science group\footnote{\url{https://wiki.helsinki.fi/display/it4sci/}} of the University of Helsinki for the high performance computing cluster Turso. TAB was partially supported by the Royal Society through the Newton International Fellowship grant n. NIF\textbackslash R1\textbackslash 201695 and by the Academy of Finland through the postdoctoral grant, decision number 330522. DL acknowledges support of NSF-DMS 1720487 and 172045. TH acknowledges support of the Emil Aaltonen Foundation junior researcher grant no. 200029 and the Academy of Finland Project 310822.

\bibliographystyle{elsarticle-num}
\bibliography{references}

\appendix
\renewcommand{\thesection}{\Alph{section}}

\section{Proof of inequality~\eqref{boundDeriv}} \label{proof.ineq}

The argument we present below follows an argument in~\cite[Prop.~10]{GL_MMNP} or \cite[Lemma 4.5]{GLintFourierOp}.

We start by observing that, for any $j \ge 0$, 
$\ell = (\ell_1,\ell_2) \in \Z^2$ with $|\ell_1|, |\ell_2| \le 2^j$, $d \in \{1,2,3\}$,  
$$\hat{\psi}_{j,\ell}^{(d)}(\xi)=W\left(2^{-2j}\xi B_{(d)}^{[\ell]}A_{(d)}^j\right)\!V_{(d)}(\xi)$$ is a smooth, bounded and compactly supported function. Hence, using the inverse Fourier transform, we can write $\psi_{j,\ell}^{(d)}(x)=\int_{\hat{\R}^4}e^{2 \pi i \xi x} \hat{\psi}_{j,\ell}^{(d)}(\xi) \, d\xi$ and, thus, we have that, for any $x \in \R^4$,  
\begin{equation}  \label{eq.est1}
   \left| \psi_{j,\ell}^{(d)}(x) \right| \leq \int_{\supp\{\hat{\psi}_{j,\ell}^{(d)}\}}\left| \hat{\psi}_{j,\ell}^{(d)}(\xi) \right| d\xi \leq m\!\left(\supp\{\hat{\psi}_{j,\ell}^{(d)}\} \right) \|\hat{\psi}_{j,\ell}^{(d)}\|_{\infty}, 
\end{equation}
where $m\!\left(\supp\{\hat{\psi}_{j,\ell}^{(d)}\}\right)$ denotes the Lebesgue measure of the support of $\hat{\psi}_{j,\ell}^{(d)}$. 

Performing integration by parts and using the regularity of $\hat{\psi}_{j,\ell}^{(d)}$, we also observe that, for any $N \in \N$,
\begin{equation}  \label{eq.est2}
    \int_{\supp\{\hat{\psi}_{j,\ell}^{(d)}\}} e^{2 \pi i \xi x} \Delta_N \left(\hat{\psi}_{j,\ell}^{(d)}\right)(\xi) \, d\xi = (2 \pi)^{2N} |x|^{2N} \psi_{j,\ell}^{(d)}(x),
\end{equation}
 where $\Delta_N=\sum_{i=1}^{4}\frac{\partial^{2N}}{\partial \xi_i^{2N}}$. Using~\eqref{eq.est1} and \eqref{eq.est2}, it follows that
 \begin{eqnarray}
 \left| \psi_{j,\ell}^{(d)}(x)\right| \left( 1+|x|^{2} \right)^N \hspace{-0.2cm} &\leq & \hspace{-0.2cm}
     \left| \psi_{j,\ell}^{(d)}(x)\right| \left( 1+ (2 \pi)^{2} |x|^{2} \right)^N \nonumber \\ 
     &\leq & \hspace{-0.2cm}
     \left| \psi_{j,\ell}^{(d)}(x)\right|  N\left( 1+ (2 \pi)^{2N} |x|^{2N} \right) \nonumber \\
     &\leq & \hspace{-0.2cm} N\, m\!\left(\supp\{\hat{\psi}_{j,\ell}^{(d)}\}\right) \left( \|\hat{\psi}_{j,\ell}^{(d)}\|_{\infty} + \|\Delta_N\hat{\psi}_{j,\ell}^{(d)}\|_{\infty} \right)\!.
     \label{eq.est3}
 \end{eqnarray}
 Using the conditions on the support of $W$ and $V_{(d)}$, we observe that for $d=1$ we have that $|\xi_2|,|\xi_3|\leq|\xi_1|$ and $2^{2j-4}\leq|2^{2j}\xi_1|, |2^{2j}\xi_4|\leq 2^{2j-1}$. Hence $|\xi_j|\leq 2^{-1}$ for $j=1,2,3,4$, which shows that $m \!\left(\supp\{\hat{\psi}_{j,\ell}^{(1)}\}\right)<C$, where $C$ is a constant independent of $j, \ell$. The same property holds for $d=2,3$ using a similar argument. Thus, combining this observation with \eqref{eq.est3}, it follows that there is a constant $C_N$, independent of $j, \ell, d$ such that
$$ \psi_{j,\ell}^{(d)}(x)\leq C_{N}(1+|x|^2)^{-N}. $$
 
To prove a similar estimate for the partial derivatives of $\psi_{j,\ell}^{(d)}$,
we start by using the properties of the Fourier transform so that, for any $\nu \in (\N \cup \{0\})^4$ and any $i \in \{1,2,3,4\}$, we  write the partial derivatives of $\psi_{j,\ell}^{(d)}$ as
\begin{eqnarray*}
    \partial_{x_i}^{\nu}\psi_{j,\ell}^{(d)}(x)=\int_{\hat{\R}^4} (2 \pi i \xi_i )^{\nu} e^{2\pi i \xi x} \hat{\psi}_{j,\ell}^{(d)}(\xi) \, d\xi.
\end{eqnarray*}
The rest of the argument is now very similar to the argument we used for $\psi_{j,\ell}^{(d)}$.

Similar to \eqref{eq.est1}, denoting $S_{j,\ell,d} =\supp\{  \hat{\psi}_{j,\ell}^{(d)}\}$, we have that 
\begin{equation*}  
   \left| \partial_{x_i}^{\nu} \psi_{j,\ell}^{(d)}(x) \right| \leq \int_{S_{j,\ell,d}}\left| (2 \pi \xi_i)^\nu \hat{\psi}_{j,\ell}^{(d)}(\xi) \right| d\xi \leq (2 \pi)^\nu m\!\left(S_{j,\ell,d} \right) \sup |\xi_i^\nu \hat{\psi}_{j,\ell}^{(d)}(\xi)|, 
\end{equation*}
where again the quantity $\sup |\xi_i^\nu \hat{\psi}_{j,\ell}^{(d)}(\xi)|$ is bounded
by a constant independently of $j,\ell,d$, due to the conditions on the support of $W$ and $V_{(d)}$.

We can similarly derive an analogues of \eqref{eq.est2} by applying integration by parts 
to the integral
$$   \int_{\supp\{\hat{\psi}_{j,\ell}^{(d)}\}} e^{2 \pi i \xi x} \Delta_N \left( (2 \pi i \xi )^{\nu} \hat{\psi}_{j,\ell}^{(d)}\right)(\xi) \, d\xi $$
and using the observation  that, for any $\beta \in (\N \cup \{0\})^4$, we have 
\begin{eqnarray*}
 \partial_{\xi_i}^{\beta} \left(  (2 \pi i \xi )^{\nu} \hat{\psi}_{j,\ell}^{(d)}(\xi)\right)= 
 \sum_{\delta+ \gamma = \beta} C_{\delta, \gamma} \, \partial_{\xi_i}^{\delta} \left( (2 \pi i \xi )^{\nu}\right) \partial_{\xi_i}^{\gamma} \left( \hat{\psi}_{j,\ell}^{(d)}(\xi)\right),
\end{eqnarray*}
where the constants $C_{\delta,\gamma}$ are independent of $j,\ell,d$. 
Finally, combining the two estimates as in \eqref{eq.est3}, we conclude that there is a constant $C_{\nu,N}$ independent of $j,\ell,d$ such that, for any $N \in \N$, we have
\begin{equation*}
    \partial_{x}^{\nu}\psi_{j,\ell}^{(d)}(x)\leq C_{\nu,N}(1+|x|^2)^{-N}.
\end{equation*}
$\qed$

\section{Proofs of Theorems~\ref{theo0} and~\ref{theo1}} \label{s.proof}

Here we assume the notation introduced in Sec.~\ref{s.main} where $f_Q = f \, w_Q$, for $f\in\E(A)$. 
We remark that the localization window $w_Q$ acts in the $x_1 x_2 x_3$ space. Below, we will select $w_Q$ 
appropriately so that we can analyze the discontinuity surface $\partial B$ locally. Recall that the surface is $C^2$ regular by hypothesis. 

Hence, by choosing $j>j_0$ sufficiently large, the scale $2^{-j}$ is small enough so that, over a cube $Q$ of side $2^{-j}$, the surface $\partial B$ may be parametrized as $x_1=E(x_2,x_3)$ or $x_2=E(x_1,x_3)$ or $x_3=E(x_1,x_2)$, where the function $E_i$, for $i=1$ or $i=2$ = $i=3$, is twice continuously differentiable.

For simplicity, we assume that this surface has parametrization  $$x_1=E(x_2,x_3), \quad |x_2|,|x_3|\leq 2^{-j},$$
as the other cases can be analyzed with a very similar argument. By a suitable translation, we may assume that the surface contains the origin, that is $k=(0,0,0)$, and the normal direction of the surface at $(0,0,0)$ is $(1,0,0)$. This is equivalent to assuming that $E(0,0)=E_{x_2}(0,0)=E_{x_3}(0,0)=0$. There is no loss in generality in analyzing only this case since the situation where the surface does not contain the origin or has a different normal direction can be easily converted into this case by translation and rotation. So, the function $f_Q$ is localized on $Q=[0,2^{-j}]^3$. To simplify notation, for a function $g(x)$ with $x\in\R^2$ and $m=(m_1,m_2)$ with $0\leq|m|=m_1+m_2\leq2$, we will write $\f{\partial^m}{\partial x^m}g$ as $g_m$.

The second order Taylor expansion of $E$ around $(0,0)$ reduces to the remainder alone, that is, $$E(x_2,x_3)=\f{1}{2}(E_{(2,0)}(c)x_2^2+2E_{(1,1)}(c)x_2x_3+E_{(0,2)}(c)x_3^2),$$ where $c=(c_2,c_3)\in[-2^{-j},2^{-j}]^2$. Therefore, for $j>j_0$ we have
$$|E(x_2,x_3)|\leq 2^{-2j}(\|E_{(2,0)}\|_\infty+\|E_{(1,1)}\|_\infty+\|E_{(0,2)}\|_\infty).$$
We will discuss the case $j\leq j_0$ further below.

Recall that $f$ has the form $f(x_1,x_2,x_3,x_4)=h(x_1,x_2,x_3)\X_B(x_1,x_2,x_3)g(x_4)$ and we want to  estimate the decay of $f$ near the surface of discontinuity.  Hence, for $j\in\Z$, we define the \textbf{surface fragment}  as the function
\begin{equation}\label{surFrag}
    h_Q(x_1,x_2,x_3)=w(2^j x) h(x_1,x_2,x_3)\X_{\{x_1>E(x_2,x_3)\}}(x_1,x_2,x_3) 
\end{equation}
where $w\in C^\infty([-1,1]^3)$ is a non-negative window function. That is,  $w(2^j\cdot)$ is the window function $w_Q$ announced above.

Note that $h\in C^2([0,1]^3)$, hence $h_Q$ is supported on $[0,2^{-j}]^3$. Consequently, we define the localized version of $f$ as $f_Q(x_1,x_2,x_3,x_4)=h_Q(x_1,x_2,x_3)g(x_4)$.

\subsection{Analysis of the Surface Fragment}
We aim at deriving $L^2$ estimates for the elements of the Parseval frame of cylindrical shearlets against $f_Q=h_Q g$, where $h_Q$ is the the surface fragment~\eqref{surFrag}. For our analysis below, it will be sufficient to consider the interior cylindrical shearlets~\eqref{shear1} associated with the pyramidal region $\P_1$. Boundary shearlets and interior shearlets in the regions $\P_2$ and $\P_3$ satisfy similar support and regularity conditions, so that the corresponding estimates against $f_Q$ are very similar.

In the following, we will express the first three coordinates of $\xi\in\R^4$ in spherical coordinates, so we write $(\xi_1,\xi_2,\xi_3)=(\rho \cos\theta \sin\phi, \rho \sin\theta \sin\phi, \rho \cos\phi) $ where $\rho>0$, $\theta\in[0,2\pi)$ and $\phi\in[0,\pi]$. Since we only consider the region $\P_1$, we can assume that $\phi\in\left[\f{\pi}{4},\f{3\pi}{4}\right]$ $\theta\in\left[-\f{\pi}{4},\f{\pi}{4}\right]$. We additionally remark that the variables $\xi_2$ and $\xi_3$ are symmetric in $\P_1$; thus, we may assume that $|\ell_1|\leq|\ell_2|$.

For $\xi\in\P_1\subset\R^4$, $j\geq0$, $|\ell_1|\leq|\ell_2|\leq2^{j}$, we let 
\begin{equation}\label{gammEq}
    \Gamma_{j,\ell}(\xi)=W(2^{-2j}\xi) \, v(2^{j}\tfrac{\xi_2}{\xi_1}-\ell_1) \, v(2^{j}\tfrac{\xi_3}{\xi_1}-\ell_2).
\end{equation}
Using this notation, the interior shearlets~\eqref{sh.one.f} associated with the pyramidal region $\P_1$ may be written as $$\hat{\psi}_{j,\ell,k}^{(1)}(\xi)=2^{-3j}\Gamma_{j,\ell}(\xi) \, e^{2\pi i\xi A_{(1)}^{-j}B_{(1)}^{[-\ell]}k}.$$

We have the following Lemma whose proof follows by direct calculation and is very similar to Sec.4.3 in~\cite{GL_3D}. Below, we use the multi-index notation $m=(m_1,m_2,m_3,m_4)\in\N^4$  with $|m|=m_1+m_2+m_3+m_4$ and write $x^m=x_1^{m_1} x_2^{m_2} x_3^{m_3} x_4^{m_4}$ and $\f{\partial^m}{\partial\xi^m}\widehat{f}=\f{\partial^{m_1}}{\partial\xi_1^{m_1}}\f{\partial^{m_2}}{\partial\xi_2^{m_2}}\f{\partial^{m_3}}{\partial\xi_3^{m_3}}\f{\partial^{m_4}}{\partial\xi_4^{m_4}}\widehat{f}$. 

\begin{lemma}\label{surfLemma}
Let 
$h_Q$ be the surface fragment defined by~\eqref{surFrag} and $\hat h_Q$ be the corresponding Fourier transform. Let 
$$(\xi_1,\xi_2,\xi_3)=(r\sin\theta'\cos\phi',-r\cos\theta',-r\sin\theta'\sin\phi') \subset \P_1$$ (which implies $|\phi'|\leq \pi/4$). We then have the following estimates.
\begin{itemize}
    \item[(a)] If the support of $h_Q$ does not intersect the surface $\partial B$, then $$\int_{2^{2j-4}}^{2j+1}\int_{0}^{2\pi}\left|\f{\partial^{m_1}}{\partial \xi_1^{m_1}}\f{\partial^{m_2}}{\partial \xi_2^{m_2}}\f{\partial^{m_3}}{\partial \xi_3^{m_3}}\widehat{h}_Q(r,\theta',\phi')\right|^2 d\theta' dr\leq C2^{-2jm_1}2^{-12j}$$
\item[(b)]     If $h_Q$ intersects the surface $\partial B$ and $|\sin\phi'|\leq 2^{1-j}$, then $$\int_{2^{2j-4}}^{2j+1}\int_{0}^{2\pi}\left|\f{\partial^{m_1}}{\partial \xi_1^{m_1}}\f{\partial^{m_2}}{\partial \xi_2^{m_2}}\f{\partial^{m_3}}{\partial \xi_3^{m_3}}\widehat{h}_Q(r,\theta',\phi')\right|^2 d\theta' dr\leq C2^{-2jm_1}2^{-7j}$$ 
\item[(c)]     If $h_Q$ intersects the surface $\partial B$ and $|\sin\phi'|\geq 2^{1-j}$, then $$\int_{2^{2j-4}}^{2j+1}\int_{0}^{2\pi}\left|\f{\partial^{m_1}}{\partial \xi_1^{m_1}}\f{\partial^{m_2}}{\partial \xi_2^{m_2}}\f{\partial^{m_3}}{\partial \xi_3^{m_3}}\widehat{h}_Q(r,\theta',\phi')\right|^2 d\theta' dr\leq C2^{-2jm_1}2^{-12j}|\sin\phi'|^{-5}$$
\end{itemize}
\end{lemma}
Note that, in the above lemma, the notation $$\f{\partial^{m_1}}{\partial \xi_1^{m_1}}\f{\partial^{m_2}}{\partial \xi_2^{m_2}}\f{\partial^{m_3}}{\partial \xi_3^{m_3}}\widehat{h}_Q(r,\theta',\phi')$$ means that we first compute $\f{\partial^{m_1}}{\partial \xi_1^{m_1}}\f{\partial^{m_2}}{\partial \xi_2^{m_2}}\f{\partial^{m_3}}{\partial \xi_3^{m_3}}\widehat{h}_Q(\xi_1,\xi_2,\xi_3)$ and next we make the change of variable $(\xi_1,\xi_2,\xi_3)=(r\sin\theta'\cos\phi',-r\cos\theta',-r\sin\theta'\sin\phi').$ 

We also have the following Lemma whose proof follows by direct calculation and is very similar to \cite[Lemma 2.5]{GL_3D}.
\begin{lemma}\label{gammaDeriv}
Let $\Gamma_{j,\ell}$ be defined by~\eqref{gammEq}. Then, for $j\geq 1$, $|\ell_1|\leq |\ell_2|\leq 2^j$ and $m=(m_1,m_2,m_3,m_4)\in\N^4$ we have 
$$\left|\f{\partial^{m_1}}{\partial \xi_1^{m_1}}\f{\partial^{m_2}}{\partial \xi_2^{m_2}}\f{\partial^{m_3}}{\partial \xi_3^{m_3}}\f{\partial^{m_4}}{\partial \xi_4^{m_4}}\Gamma_{j,\ell}(\xi)\right|\leq C_m 2^{-m_1 j}2^{-|m|j}(1+|\ell_2|)^{m_1}.$$
\end{lemma}

Observing that the supports of $\Gamma_{j,\ell_1,\ell_2}$ and $\Gamma_{j,\ell_1',\ell_2}$ are disjoint provided that $\ell_1\neq\ell_1'$, it follows from Lemma~\ref{gammaDeriv} that $$\sum_{\ell_1=-|\ell_2|}^{|\ell_2|}\left|\f{\partial^{m_1}}{\partial \xi_1^{m_1}}\f{\partial^{m_2}}{\partial \xi_2^{m_2}}\f{\partial^{m_3}}{\partial \xi_3^{m_3}}\f{\partial^{m_4}}{\partial \xi_4^{m_4}}\Gamma_{j,\ell}(\xi)\right|\leq C_m 2^{-m_1 j}2^{-|m|j}(1+|\ell_2|)^{m_1}.$$

We can now prove the following result.

\begin{theorem}\label{sumDerivBoundTheo}
 Let $f_Q=h_Q g$, where $h_Q$ is the surface fragment given by~\eqref{surFrag}, $g\in C^2([-1,1])$ and $\Gamma_{j,\ell}$ is given by~\eqref{gammEq}. Let  $m_f=(m_{f_1},m_{f_2},m_{f_3},m_{f_4})$ and $m_\gamma=(m_{\gamma_1},m_{\gamma_2},m_{\gamma_3},m_{\gamma_4})$ be multi-indexes. Then there exists a constant $C$ independent of $j, \ell$ such that
\begin{eqnarray*}
  &&\sum_{\ell_1=-|\ell_2|}^{|\ell_2|}\int_{\R^4} \left| \f{\partial^{m_{f}}}{\partial\xi^{m_f}}\widehat{f}_Q(\xi)\right|^2 \left|\f{\partial^{m_{\gamma}}}{\partial\xi^{m_\gamma}}\Gamma_{j,\ell}(\xi)\right|^2 d\xi\\
  &&\leq C  2^{-m_{\gamma_1} j}2^{-|m_\gamma|j}(1+|\ell_2|)^{m_{\gamma_1}} 2^{-2jm_{f_1}}\left(2^{-4j}\left(1+ |\ell_2|^{-5}\right) +2^{-9j}\right).
 \end{eqnarray*}
\end{theorem}
\textbf{Proof.} 
Recall that the support of $\Gamma_{j,\ell}$ is contained in $\P_1$ and depends on the supports of $W$ and $v$. So, for  $\xi=(\xi_1,\xi_2,\xi_3,\xi_4)\in\supp\Gamma_{j,\ell}$, we have $\xi_i\in\left[-2^{2j-1},2^{2j-1}\right]\setminus\left[-2^{2j-4},2^{2j-4}\right]$ for $i=1,2,3,4$,  $\left|2^j\f{\xi_2}{\xi_1}-\ell_1\right|\leq 1$ and $\left|2^j\f{\xi_3}{\xi_1}-\ell_2\right|\leq 1$. By applying a change of variables into spherical coordinates, we can write $(\xi_1,\xi_2,\xi_3)=(r\sin\theta'\cos\phi',-r\cos\theta',-r\sin\theta'\sin\phi')$, where we have $\left|-2^j\f{\cot\theta'}{\cos\phi'}-\ell_1\right|\leq 1$ and $\left|-2^j\tan\phi'-\ell_2\right|\leq 1$. Thus, $$r^2=\xi_1^2+\xi_2^3+\xi_3^2=\xi_1^2\left(1+(\tfrac{\cot\theta'}{\cos\phi'})^2+\left(\tan\phi'\right)^2\right)$$ and $2^{2j-4}\leq r\leq 2^{2j+2}$. We also remark that $|\phi'|\leq\pi/4$ since $\Gamma_{j,\ell}$ is supported on $\P_1$. In addition, from $(-1-\ell_2)2^{-j}\leq \tan\phi'\leq (1-\ell_2)2^{-j}$ and the Taylor expansion of the tangent function, we see that $\phi'$ must be contained in an interval $I_{\phi'}$ of length $C2^{-j}$. Hence, using Lemma~\ref{gammaDeriv} we have
\begin{eqnarray*}
    && \sum_{\ell_1=-|\ell_2|}^{|\ell_2|}\int_{\widehat{\R}^4}\left| \f{\partial^{m_{f}}}{\partial\xi^{m_f}}\widehat{f}_Q(\xi)\right|^2 \left|\f{\partial^{m_{\gamma}}}{\partial\xi^{m_\gamma}}\Gamma_{j,\ell}(\xi)\right|^2 d\xi \\
    &\leq& C_{m_\gamma} 2^{-m_{\gamma 1} j}2^{-|m_\gamma|j}(1+\ell_2)^{m_{\gamma 1}}\int_{\widehat{\R}^4}\left| \f{\partial^{m_{f}}}{\partial\xi^{m_f}}\widehat{f}_Q(\xi)\right|^2 d\xi\\
    & \leq& C_{m_\gamma} 2^{-m_{\gamma_1} j}2^{-|m_\gamma|j}(1+\ell_2)^{m_{\gamma_1}}
    \int_{\widehat{\R}^4}\left| \f{\partial^{m_{f_1} }}{\partial \xi_1^{m_{f_1}}}\f{\partial^{m_{f_2}}}{\partial \xi_2^{m_{f_2}}}\f{\partial^{m_{f_3}}}{\partial \xi_3^{m_{f_3}}}\widehat{h}_Q(\xi_1,\xi_2,\xi_3) \right|^2 \!\!d\xi\\
    &=& C_{m_\gamma} 2^{-m_{\gamma_1} j}2^{-|m_\gamma|j}(1+\ell_2)^{m_{\gamma_1}}
    \int_{I_{\phi'}}\int_{2^{2j-4}}^{2^{2j+2}} \int_{0}^{2\pi} r^2 |\sin\theta'|\\
    &\times& \left|\f{\partial^{m_{f_1} }}{\partial \xi_1^{m_{f_1}}}\f{\partial^{m_{f 2}}}{\partial \xi_2^{m_{f_2}}}\f{\partial^{m_{f 3}}}{\partial \xi_3^{m_{f_3}}}\widehat{h}_Q(r,\theta',\phi') \right|^2   d\theta' dr d\phi'\\
    & \leq&  C_{m_\gamma} 2^{-m_{\gamma_1} j}2^{-|m_\gamma|j}(1+\ell_2)^{m_{\gamma_1}} 2^{4j} \\ 
    & \times &\int_{I_{\phi'}}\int_{2^{2j-4}}^{2^{2j+2}} \!\!\int_{0}^{2\pi}\left| \f{\partial^{m_{f_1} }}{\partial \xi_1^{m_{f_1}}}\f{\partial^{m_{f_2}}}{\partial \xi_2^{m_{f_2}}}\f{\partial^{m_{f_3}}}{\partial \xi_3^{m_{f_3}}}\widehat{h}_Q(r,\theta',\phi') \right|^2  d\theta' dr d\phi'.
\end{eqnarray*}
Next, we apply Lemma~\ref{surfLemma}. In the  no-intersection case we have
\begin{eqnarray*}
    && \sum_{\ell_1=-|\ell_2|}^{|\ell_2|}\int_{\widehat{\R}^4}\left| \f{\partial^{m_{f}}}{\partial\xi^{m_f}}\widehat{f}_Q(\xi)\right|^2 \left|\f{\partial^{m_{\gamma}}}{\partial\xi^{m_\gamma}}\Gamma_{j,\ell}(\xi)\right|^2  d\xi \\
    &&  \leq   C_{m_\gamma} 2^{-m_{\gamma_1} j}2^{-|m_\gamma|j}(1+|\ell_2|)^{m_{\gamma_1}} 2^{4j}\int_{I_{\phi'}} 2^{-2j m_{f_1}}2^{-12j}d\phi'\\
    &&= C_{m_\gamma} 2^{-m_{\gamma_1} j}2^{-|m_\gamma|j}(1+|\ell_2|)^{m_{\gamma_1}} 2^{-2jm_{f_1}}2^{-9j}.
\end{eqnarray*}
In the intersection case,  if $|\sin\phi'|\leq 2^{1-j}$, we have 
\begin{eqnarray*}
&& \sum_{\ell_1=-|\ell_2|}^{|\ell_2|}\int_{\widehat{\R}^4}\left| \f{\partial^{m_{f}}}{\partial\xi^{m_f}}\widehat{f}_Q(\xi)\right|^2 \left|\f{\partial^{m_{\gamma}}}{\partial\xi^{m_\gamma}}\Gamma_{j,\ell}(\xi)\right|^2 d\xi \\
&&\leq C_{m_\gamma} 2^{-m_{\gamma_1} j}2^{-|m_\gamma|j}(1+|\ell_2|)^{m_{\gamma_1}} 2^{4j} \int_{I_{\phi'}}2^{-2jm_1}2^{-7j}d\phi'\\
    &&= C_{m_\gamma} 2^{-m_{\gamma_1} j}2^{-|m_\gamma|j}(1+|\ell_2|)^{m_{\gamma_1}} 2^{-2jm_{f_1}}2^{-4j}.
\end{eqnarray*}    
On the other hand, if $|\sin\phi'|\geq 2^{1-j}$ (in which case $2^j|\sin\phi'|$ is equivalent to $|\ell_2|$), then we have
\begin{eqnarray*}
&& \sum_{\ell_1=-|\ell_2|}^{|\ell_2|}\int_{\widehat{\R}^4}\left| \f{\partial^{m_{f}}}{\partial\xi^{m_f}}\widehat{f}_Q(\xi)\right|^2 \left|\f{\partial^{m_{\gamma}}}{\partial\xi^{m_\gamma}}\Gamma_{j,\ell}(\xi)\right|^2 d\xi \\ 
&&\leq C_{m_\gamma} 2^{-m_{\gamma_1} j}2^{-|m_\gamma|j}(1+|\ell_2|)^{m_{\gamma_1}} 2^{4j}
    \int_{I_{\phi'}}2^{-2jm_1}2^{-12j}|\sin\phi'|^{-5}d\phi'\\
    &&\leq C_{m_\gamma} 2^{-m_{\gamma_1} j}2^{-|m_\gamma|j}(1+|\ell_2|)^{m_{\gamma_1}} 2^{-2jm_{f_1}}2^{-4j} |\ell_2|^{-5}.
\end{eqnarray*}
This proves the theorem.$\qed$

To prove Theorem~\ref{theo0} we modify an idea from~\cite{GL_3D} to take into account the fourth variable. For that, we introduce the following differential operator:
\begin{equation}\label{diffOperator}
   L=\left(I-(\tfrac{2^{2j}}{2\pi (1+|\ell_2|)})^2 \tfrac{\partial^2}{\partial \xi_1^2}\right)\!\! \left( I- (\tfrac{2^j}{2\pi})^2 \tfrac{\partial^2}{\partial \xi_2^2}\right) \!\!\left( I- (\tfrac{2^j}{2\pi})^2 \tfrac{\partial^2}{\partial \xi_3^2}\right) \!\! \left(I-(\tfrac{2^{2j}}{2 \pi})^2\tfrac{\partial^2}{\partial\xi_4^2}\right) 
\end{equation}

Using Theorem~\ref{sumDerivBoundTheo}, a direct computation gives the following result.
\begin{theorem}\label{L.inequatily}
 Let $f=h_Q \, g$, where $h_Q$ is a surface fragment given by~\eqref{surFrag}, and $\Gamma_{j,\ell}$ be given by~\eqref{gammEq}. Then, for $j\geq 0$ and $|\ell_2|\leq 2^j$, we have
 $$\sum_{\ell_1=-|\ell_2|}^{|\ell_2|} \int_{\widehat{\R}^4}\left|L\left(\widehat{f}(\xi)\Gamma_{j,\ell}(\xi)\right)\right|^2 d\xi\leq C 2^{-4j}(1+|\ell_2|)^{-5}.$$
\end{theorem}
As observed above, Theorem~\ref{L.inequatily} gives an estimate valid for $f$ in the region $\P_1$. A very similar estimate can be derived for the regions $\P_2$ and $\P_3$, using appropriate modifications of the differential operator $L$. 

\subsection{Proof of Theorem~\ref{theo0}}
Fix $j\geq0$. By our remark above, it is enough to consider the region $\P_1$ only. For $\mu\in M_j$,  the shearlet coefficients of $f_Q$ associated to $\P_1$ can be written as 
$$\ip{f_Q}{\tilde{\psi}_\mu}=\ip{f_Q}{\psi_{j,k,\ell}^{(1)}}=|\det A_{(1)}|^{-j/2}\int_{\widehat{R}^4}\widehat{f_Q}(\xi)\Gamma_{j,\ell}(\xi)e^{2\pi i \xi A_{(1)}^{-j}B_{(1)}^{-[\ell]}k}d\xi$$
where $\Gamma_{j,\ell}$ is given by~\eqref{gammEq}. Using the equivalent definition of the weak $\ell^1$ norm, we need to show 
\begin{eqnarray}\label{ell1Ineq}
    \#\{\mu\in M_j\,:\,|\ip{f_Q}{\tilde{\psi}_\mu}|>\epsilon\}\leq C 2^{-2j}\epsilon^{-1}.
\end{eqnarray}
We observe that $\xi A_{(1)}^{-j}B_{(1)}^{-[\ell]}k=(k_1-k_2\ell_1-k_3\ell_2)2^{-2j}\xi_1+k_2 2^{-j}\xi_2+k_3 2^{-j}\xi_3+ k_4 2^{-2j} \xi_4$. Hence, letting $L$ to be the differential operator in~\eqref{diffOperator}, we have that  
$$ L\left( e^{2\pi i \xi A_{(1)}^{-j}B_{(1)}^{-[\ell]}k} \right) = \begin{cases}
  G_0(k,\ell) \, e^{2\pi i \xi A_{(1)}^{-j}B_{(1)}^{-[\ell]}k} & \text{if } \ell_2 = 0 \\
  G_1(k,\ell) \, e^{2\pi i \xi A_{(1)}^{-j}B_{(1)}^{-[\ell]}k}& \text{if } \ell_2\neq 0
\end{cases}
$$
where
\begin{eqnarray*}
   G_0(k,\ell) &=& G_0(k) =  (1+k_1^2)(1+k_2^2) (1+k_3)^2 (1+k_4)^2\\
   G_1(k, \ell) &=&  (1+(\tfrac{|\ell_2|}{1+|\ell_2|})^2
(\tfrac{k_1}{|\ell_2|}-\tfrac{k_2 \ell_1}{|\ell_2|}\pm k_3)^2) \, (1+k_2^2) (1+k_3)^2 (1+k_4)^2.
\end{eqnarray*}
The $\pm$ sign in the above expression follows from dividing $\ell_2$ by $|\ell_2|$; in other words $\pm k_3 = \mbox{sign}(\ell_2) k_3$.

Hence, a direct computation (using integration by parts) shows that 
\begin{equation}  \label{eq.ipg}
\ip{f_Q}{\tilde{\psi}_\mu}=|\det(A_{(1)}|^{-\frac j2}\int_{\widehat{\R}^4}L\left(\widehat{f_Q}(\xi)\Gamma_{j,\ell}(\xi)\right) G_i(\ell, k)^{-1}
 e^{2\pi i \xi A_{(1)}^{-j}B_{(1)}^{-[\ell]}k} \, d\xi,
 \end{equation}
where  $G_i=G_0$ if $\ell_2 = 0$ and $G_i=G_1$ if $\ell_2 \neq 0$.

We next consider the cases $\ell_2=0$ and $\ell_2\neq0$ separately. 

For $\ell_2\neq0$, \eqref{eq.ipg} gives
$$G_1(k,\ell)\ip{f_Q}{\tilde{\psi}_\mu}=|\det(A_{(1)}|^{-j/2}\int_{\widehat{\R}^4}L\left(\widehat{f_Q}(\xi)\Gamma_{j,\ell}(\xi)\right) e^{2\pi i \xi A_{(1)}^{-j}B_{(1)}^{-[\ell]}k} d\xi.$$
Let $K=(K_1, K_2 ,K_3, K_4)\in\Z^4$ and define
\begin{eqnarray*}
    R_K &=& \{k=(k_1,k_2,k_3,k_4)\in\Z^4:\,  \f{k_1}{|\ell_2|}\in[K_1,K_1+1],\, -\f{k_2 \ell_1}{|\ell_2|}\in[K_2,K_2+1], \\
    && k_3=K_3,\, k_4=K_4\}.
\end{eqnarray*}
For fixed $j,$ $\ell$, the set $\{|\det A_{(1)}|^{-j/2} e^{2\pi i \xi A_{(1)}^{-j}B_{(1)}^{-[\ell]}k}\!: k\in\Z^4\}$ is an orthonormal basis for $L^2$ functions defined on $[-1/2,1/2]^4 B_{(1)}^{[\ell]}A_{(1)}^{j}$ (which contains the support of $\Gamma_{j,\ell}$). It follows that
\begin{eqnarray*}
    && \sum_{k\in R_K} \!\! G_1(k,\ell)^2|\ip{f_Q}{\tilde{\psi}}|^2 \\
    && =  \sum_{k\in R_K}|\det(A_{(1)})|^{-j}\left|\int_{\widehat{\R}^4}L\left(\widehat{f_Q}(\xi)\Gamma_{j,\ell}(\xi)\right) e^{2\pi i \xi A_{(1)}^{-j}B_{(1)}^{-[\ell]}k} d\xi\right|^2\\
    && \leq \|L(\widehat{f_Q}\Gamma_{j,\ell})\|^2\\
    && = \int_{\widehat{\R}^4} \left| L\left(\widehat{f_Q}(\xi)\Gamma_{j,\ell}(\xi)\right)\right|^2 d\xi.
\end{eqnarray*}
Therefore, 
\begin{eqnarray*}
   && \!\! \sum_{\ell_1=-|\ell_2|}^{|\ell_2|}\sum_{k\in R_K}|\ip{f_Q}{\tilde{\psi}}|^2 \\&=& \!\! \sum_{\ell_1=-|\ell_2|}^{|\ell_2|}\sum_{k\in R_K} \!\! G_1(k,\ell)^{-2}|\det(A_{(1)}|^{-j}\left|\int_{\widehat{\R}^4} \!\! L\left(\widehat{f_Q}(\xi)\Gamma_{j,\ell}(\xi)\right) e^{2\pi i \xi A_{(1)}^{-j}B_{(1)}^{-[\ell]}k} d\xi\right|^2\\
    &\leq& C\left(1+(K_1+K_2\pm K_3)^2\right)^{-2} (1+K_2^2)^{-2} (1+K_3^2)^{-2} (1+K_4^2)^{-2} \\ &\times & \sum_{\ell_1=-|\ell_2|}^{|\ell_2|}\int_{\widehat{\R}^4} \left|  L\left(\widehat{f_Q}(\xi)\Gamma_{j,\ell}(\xi)\right)\right|^2 d\xi.
\end{eqnarray*}
Using Theorem~\ref{L.inequatily}, we have 
\begin{eqnarray*}
    \sum_{\ell_1=-|\ell_2|}^{|\ell_2|}\sum_{k\in R_K}|\ip{f_Q}{\tilde{\psi}}|^2\leq H_K^{-2} C 2^{-4j}(1+|\ell_2|)^{-5},
\end{eqnarray*}
where we define $H_K= \left(1+(K_1+K_2\pm K_3)^2\right) (1+K_2^2) (1+K_3^2) (1+K_4^2)$. 

For fixed $j$ and $\ell$, we let $R_{K,\epsilon}=\{k\in R_K:\, |\ip{f}{\psi_{j,k,\ell}^{(1)}}|>\epsilon\}$ and $N_{j,\ell,K}(\epsilon)=\# R_{K,\epsilon}$. Then, by the condition $|\ell_1|\leq |\ell_2|$, we have that $N_{j,k,\ell}(\epsilon)\leq C(1+|\ell_2|)^2$ and, thus, $\sum_{\ell_1=-|\ell_2|}^{|\ell_2|}N_{j,k,\ell}(\epsilon)\leq (1+|\ell_2|)^3$. Now, for $k\in R_K$ such that $|\ip{f_Q}{\psi_{j,k,\ell}^{(1)}}|>\epsilon$, we have
\begin{eqnarray*}
    \epsilon^2 N_{j,\ell,K}(\epsilon) \leq \sum_{k\in R_{K,\epsilon}}|\ip{f_Q}{\tilde{\psi}}|^2\leq \sum_{k\in R_K}\leq |\ip{f_Q}{\tilde{\psi}}|^2,
\end{eqnarray*}
which implies
\begin{eqnarray*}
    \sum_{\ell_1=-|\ell_2|}^{|\ell_2|}N_{j,\ell,K}(\epsilon)\leq C \, H_K^{-2} \,  2^{-4j}(1+|\ell_2|)^{-5}\epsilon^{-2}.
\end{eqnarray*}
Hence
\begin{eqnarray*}
\sum_{\ell_1=-|\ell_2|}^{|\ell_2|}N_{j,\ell,K}(\epsilon)\leq C \min\left((1+|\ell_2|)^3, H_K^{-2}  2^{-4j}(1+|\ell_2|)^{-5}\epsilon^{-2}\right).
\end{eqnarray*}
Now, let $\ell_2^*$ be defined by   $(\ell_2^*+1)^3= H_K^{-2} 2^{-4j} \epsilon^{-2} (1+\ell_2^*)^{-5}$, so $(1+\ell_2^*)^{4}=  H_K^{-1}2^{-2j}\epsilon^{-1}$. Then,
\begin{eqnarray*}
&& \sum_{\ell_2=-2^j}^{2^j}\sum_{\ell_1=-|\ell_2|}^{|\ell_2|}N_{j,\ell,K}(\epsilon) \\
&&\leq  \sum_{|\ell_2|\leq(\ell_2^*+1)}\sum_{\ell_1=-|\ell_2|}^{|\ell_2|}N_{j,\ell,K}(\epsilon)+\sum_{|\ell_2|>(\ell_2^*+1)}\sum_{\ell_1=-|\ell_2|}^{|\ell_2|}N_{j,\ell,K}(\epsilon)\\
&& \leq C\sum_{|\ell_2|\leq(\ell_2^*+1)}(|\ell_2|+1)^3+C\sum_{|\ell_2|>(\ell_2^*+1)}H_K^{-2}2^{-4j}\epsilon^{-2} (1+|\ell_2|)^{-5}\\
&&\leq C(\ell_2^*+1)^4+C H_K^{-2}2^{-4j}\epsilon^{-2}(1+\ell_2^*)^{-4}\\
&&\leq C H_K^{-1}2^{-2j}\epsilon^{-1}.
\end{eqnarray*}
Notice also that $\sum_{K\in\Z^4}H_K^{-1}<\infty$. Thus
\begin{eqnarray*}
    \#\{\mu\in M_j:\, |\ip{f_Q}{\tilde{\psi}_\mu}|>\epsilon\}&\leq&\sum_{K\in\Z^4}\sum_{\ell_2=-2^j}^{2^j}\sum_{\ell_1=-|\ell_2|}^{|\ell_2|}N_{j,\ell,K}(\epsilon)\\&\leq & C 2^{-2j}\epsilon^{-1} \sum_{K\in\Z^4}H_K^{-1}\leq C 2^{-2j}\epsilon^{-1}
\end{eqnarray*}
which gives~\eqref{ell1Ineq}. 

Next we consider the case $\ell_2=0$. In this case, 
\eqref{eq.ipg} gives
$$G_0(k)\ip{f_Q}{\tilde{\psi}_\mu}=|\det(A_{(1)}|^{-j/2}\int_{\widehat{\R}^4}L\left(\widehat{f_Q}(\xi)\Gamma_{j,\ell}(\xi)\right) e^{2\pi i \xi A_{(1)}^{-j}B_{(1)}^{-[\ell]}k} d\xi.$$
Using Theorem~\ref{L.inequatily}, we have
 that \begin{eqnarray*}
    \sum_{k\in\Z^4} G_0(k)^{2} \, |\ip{f_Q}{\psi_{j,\ell,k}^{(1)}}|^2=\int_{\widehat{\R}^4}\left|L\left(\widehat{f_Q}(\xi) \Gamma_{j,\ell}(\xi)\right)\right|^2d\xi\leq C \, 2^{-4j}.
\end{eqnarray*}
In particular, for each $k\in\Z^4$ we have $|\ip{f_Q}{\psi_{j,\ell,k}^{(1)}}|\leq C \, G_0(k)^{-1} 2^{-2j}$. Since  $\sum_{k\in\Z^4} G_0(k)^{-1}<\infty$, there is a constant $C$ such that 
\begin{eqnarray*}
\|\ip{f_Q}{\psi_{j,\ell,k}^{(1)}}\|_1=\sum_{k\in\Z^4}|\ip{f_Q}{\psi_{j,\ell,k}^{(1)}}|\leq C \, 2^{-2j}
\end{eqnarray*}
which implies $\|\ip{f_Q}{\psi_{j,\ell,k}^{(1)}}\|_{w\ell^1}\leq C \, 2^{-2j}$. This concludes the proof of the theorem when $j \leq j_0$. 

\subsection{Analysis of the coarse scale.}
At the beginning Section~\ref{s.proof}, we assumed $j>j_0$ for some $j_0>0$. Here we consider the {\it coarse scale} case $j\leq j_0$. 

We recall that $f_Q(x)=h_Q(x_1,x_2,x_3) \, g(x_4)$ where $$h_Q(x_1,x_2,x_3)=w(2^j x) \, h(x_1,x_2,x_3) \, \X_{\{x_1>E(x_2,x_3)\}}(x_1,x_2,x_3)$$ with $w\in C^\infty([-1,1]^3)$, $h\in C^2([0,1]^3)$, $g\in C^2([0,1])$. Therefore, observing that $\supp h_Q\in [-2^{-j},2^{-j}]^3$ and $g$ is also compactly supported, we have 
$$\|f_Q\|_2^2\leq \int_{\supp(h_Q)\times \supp(g)}|f_Q(x)|^2 dx\leq C \, 2^{-3j}.$$
The last inequality implies that $\|\ip{f_Q}{\tilde{\psi}_\mu}\|_{\ell^2}\leq \|f_Q\|_2\leq C \, 2^{-3j/2}$. We also notice that $$\|\ip{f_Q}{\tilde{\psi}_\mu}\|_{\ell^p}\leq N^{1/p-1/2}\|\ip{f_Q}{\tilde{\psi}_\mu}\|_{\ell^2}$$
is valid for any sequence $\{\ip{f_Q}{\tilde{\psi}_\mu}\}$ of $N$ elements. Since, at scale $2^{-j}$, there are about $2^{2j}$ shearlet elements in $Q_j^0$, so we conclude that there is a constant $C$ independent of $Q$ and $j$ such that  $$\|\ip{f_Q}{\tilde{\psi}_\mu}\|_{\ell^1}\leq C \, 2^{2j(1-1/2)}2^{-3/2j}=C \, 2^{-j/2}.$$

This completes the proof of Theorem~\ref{theo0} for $j>j_0$.

\subsection{Proof of Theorem~\ref{theo1}.}

We again write $h_Q(x_1,x_2,x_3)=h(x_1,x_2,x_3) w_Q(x_1,x_2,x_3)$ where we now assume $Q\in\Q_j^1$. With this notation, we write the localized function $f_Q$ as $f_Q(x_1,x_2,x_3,x_4)= h_Q(x_1,x_2,x_3) g(x_4)$. The following two lemmata can be proved 
using an argument very similar to Lemma 4.8 and Lemma 4.9 in~\cite{GL_3D}.
\begin{lemma}\label{fFourierBound}
Let $f_Q=f w_Q$ where $f=h\X_{B}g\in\E(A)$ is given by~\eqref{cartLike}, $Q\in \Q_j^1$ and $U_{j,\ell}$ be given by~\eqref{U.j.ell}. Then, 
\begin{eqnarray*}
    \int_{U_{j,\ell}}|\widehat{f}_Q(\xi)|^2d\xi\leq C 2^{-11j}.
\end{eqnarray*}
\end{lemma}

\begin{lemma}\label{gammaDerivLemma}
Let $m=(m_1,m_2,m_3,m_4)\in\N^4$, $\xi=(\xi_1,\xi_2,\xi_3,\xi_4)\in \R^4$ and $\Gamma_{j,\ell}$ given by~\eqref{gammEq} where $\ell=(\ell_1,\ell_2)$. Then, $$\sum_{\ell_1=2^{j}}^{2^{j}}\sum_{\ell_2=-2^{j}}^{2^{j}}\left|\f{\partial^m}{\partial\xi^m}\Gamma_{j\ell}(\xi)\right|^2\leq C_m 2^{-|m|j},$$
where $C_m$ is independent of $j$ and $\xi$, and $|m|=m_1+m_2+m_3+m_4$.
\end{lemma} 
By Lemmata~\ref{fFourierBound} and~\ref{gammaDerivLemma},  using an argument similar to the proof of Lemma 4.10 in~\cite{GL_3D}, we have the following result.
\begin{lemma}\label{T.lemma}
Let $f_Q=f w_Q$ where $f\in\E^2(A)$ and $Q\in\Q_j^1$ and set \begin{eqnarray}\label{deltaOper}
    T=\left( I-\f{2^{2j}}{(2\pi)^2}\Delta\right)
\end{eqnarray}
where $\Delta=\f{\partial^2}{\partial\xi_1^2}+\f{\partial^2}{\partial\xi_2^2}+\f{\partial^2}{\partial\xi_3^2}+\f{\partial^2}{\partial\xi_4^2}$. Then, 
\begin{eqnarray*}
    \int_{\widehat{\R}^4} \sum_{\ell_1=2^{j}}^{2^{j}}\sum_{\ell_2=-2^{j}}^{2^{j}}\left|T^2\left(\widehat{f}_Q\Gamma_{j,\ell}\right)(\xi)\right|^2 d\xi \leq C 2^{-11j}.
    \end{eqnarray*}
\end{lemma}

Now we prove Theorem~\ref{theo1}. 

\textbf{Proof of Theorem~\ref{theo1}.} 

As observed above, it will be sufficient to consider the system of interior shearlets in the pyramidal region $\P_1$ as the other pyramidal regions and the boundary shearlets can be handled in a similar way.

For $T$ given by~\eqref{deltaOper}, denoting $E_1(\xi) = e^{2\pi i \xi A_{(1)}^{-j}B_{(1)}^{-[\ell]}k}$, we have 
\begin{eqnarray*}
&& T\left(E_1(\xi)\right) = \left(1+2^{-2j}(k_1-\ell_1 k_2-\ell_2 k_3)^2+k_2^2+k_3^2+2^{-2j}k_4^2\right)E_1(\xi) \\
&& T^2\left(E_1(\xi)\right)=\left(1+2^{-2j}(k_1-\ell_1 k_2-\ell_2 k_3)^2+k_2^2+k_3^2+2^{-2j}k_4^2\right)^2 E_1(\xi).
\end{eqnarray*}
For a fixed $j\geq 0$ and $f_Q=f w_Q$ where $f=h g$ and $Q\in\Q_j^1$, using integration by parts we have 
\begin{eqnarray*}
\ip{f_Q}{\tilde{\psi}_\mu} \!\!&=&|\det A_{(1)}|^{-j/2}\int_{\widehat{\R}^4}\widehat{f}_Q(\xi)\Gamma_{j,\ell}(\xi) e^{2\pi i \xi A_{(1)}^{-j}B_{(1)}^{-[\ell]}k}d\xi\\
&=&|\det A_{(1)}|^{-j/2} (1+2^{-2j}(k_1-\ell_1 k_2-\ell_2 k_3)^2+k_2^2+k_3^2+2^{-2j}k_4^2)^{-2}\\
&\times&\int_{\widehat{\R}^4} T^2\left(\widehat{f}_Q(\xi)\Gamma_{j,\ell}(\xi)\right) e^{2\pi i \xi A_{(1)}^{-j}B_{(1)}^{-[\ell]}k}d\xi.
\end{eqnarray*}
Now, for $K=(K_1,K_2,K_3,K_4)\in\Z^4$ we set 
\begin{eqnarray*}
R_K&=&\{(k_1,k_2,k_3,k_4)\in\Z^4:\, k_3=K_3,\, k_2=K_2,\\
&& 2^{-j} k_4 \in[K_4,K_4+1],\, 2^{-j}(k_1-K_2\ell_1-K_3\ell_2)\in[K_1,K_1+1]\}.
\end{eqnarray*}
We observe that $2^j K_1\leq k_1-K_2\ell_1-K_3\ell_2\leq 2^{j}(K_1+1)$ and $2^j K_4\leq k_4\leq 2^j(K_4+1)$, so for each $K$ and $\ell$ there are only $1+2^j$ choices for $k_1$ and $k_4$ in $R_K$. Thus, the number of elements of $R_K$ is bounded by $(2^j+1)^2$. We next use an argument similar to the proof of Theorem~\ref{theo0} above. We observe that, for fixed $j$ and $\ell$ the set $\{|\det A_{(1)}|^{-j/2} e^{2\pi i \xi A_{(1)}^{-j}B_{(1)}^{-[\ell]}k}:k\in\Z^4 \}$ is an orthonormal basis for the $L^2$ functions supported on $[-1/2,1/2]B_{(1)}^{[\ell]}A_{(1)}^{j}$. Furthermore, we note that $\Gamma_{j,\ell}$ is supported on $[-1/2,1/2]B_{(1)}^{[\ell]}A_{(1)}^{j}$. Thus
\begin{eqnarray*}
    \sum_{k\in R_K}\!\!|\ip{f_Q}{\tilde{\psi}_\mu}|^2 \hspace{-0.2cm}&\!=\!&\hspace{-0.2cm} 2^{-6j}\!\!\sum_{k\in R_K}\!\!\left(1+2^{-2j}(k_1-\ell_1 k_2-\ell_2 k_3)^2+k_2^2+k_3^2+2^{-2j}k_4^2\right)^{-2}\\
    &\times& \hspace{-0.2cm}\left|\int_{\widehat{\R}^4} T^2\left(\widehat{f}_Q(\xi)\Gamma_{j,\ell}(\xi)\right) e^{2\pi i \xi A_{(1)}^{-j}B_{(1)}^{-[\ell]}k}d\xi\right|^2\\
    &\leq& \hspace{-0.2cm} 2^{-6j}(1+K_1^2+K_2^2+K_3^2+K_4^2)^{-4}\int_{\widehat{\R}^4}\left|T^2\left(\widehat{f}_Q(\xi)\Gamma_{j,\ell}(\xi)\right) \right|^2\!d\xi.
\end{eqnarray*}
From the last inequality, using lemma~\ref{T.lemma} we see that
\begin{eqnarray}
    \sum_{\ell_1, \ell_2=-2^j}^{2^j} \sum_{k\in R_K}|\ip{f_Q}{\tilde{\psi}_\mu}|^2 &\leq&C (1+K_1^2+K_2^2+K_3^2+K_4^2)^{-4} \nonumber \\
    &\times&\int_{\widehat{\R}^4}\sum_{\ell_1, \ell_2=-2^j}^{2^j}\left|T^2\left(\widehat{f}_Q(\xi)\Gamma_{j,\ell}(\xi)\right) \right|^2d\xi  \nonumber\\
    &\leq&C (1+K_1^2+K_2^2+K_3^2+K_4^2)^{-4} 2^{-11j}. \label{eq.l2est}
\end{eqnarray}
By the H\"older inequality, we have that for any $N\in\N$ \begin{eqnarray*}
    \sum_{m=1}^N|a_m|\leq\left(\sum_{m=1}^N |a_m|^2\right)^{1/2} N^{1/2}.
\end{eqnarray*}
Thus, using the last inequality with \eqref{eq.l2est} and the observation that the number of elements of $R_K$ is bounded by $(1+2^j)^2$, we have
\begin{eqnarray*}
\sum_{\ell_1, \ell_2=-2^j}^{2^j} \sum_{k\in R_K}|\ip{f_Q}{\tilde{\psi}_\mu}|&\leq& \left(2^{4j}\right)^{1/2}\left(\sum_{\ell_1, \ell_2=-2^j}^{2^j} \sum_{k\in R_K}|\ip{f_Q}{\tilde{\psi}_\mu}|^2\right)^{1/2}\\
&\leq&2^{2j} C (1+K_1^2+K_2^2+K_3^2+K_4^2)^{-2} 2^{-11j/2}.
\end{eqnarray*}
This shows that, for $f_Q=f w_Q$ with $Q\in\Q_j^1$ we have $\sum_{\mu\in M_j}|\ip{f_Q}{\tilde{\psi}_\mu}|\leq C 2^{-7j/2}$. This completes the proof of theorem. $\qed$

\end{document}